\documentclass[a4paper,12pt,reqno]{amsart}
\usepackage[utf8]{inputenc}
\usepackage[finnish,english]{babel}
\usepackage{amsthm,amssymb,mathtools,hyperref}
\usepackage{amsfonts,mathrsfs}

\mathtoolsset{showonlyrefs}

\title[Geometric inverse problems on gas giants]{Geometric inverse problems on gas giants}
\author[de Hoop]{Maarten V. de Hoop}
\address{Simons Chair in Computational and Applied Mathematics and Earth Science, Rice University}
\author[Ilmavirta]{Joonas Ilmavirta}
\address{Department of Mathematics and Statistics, University of Jyv\"askyl\"a}
\author[Kykk\"anen]{Antti Kykk\"anen}
\address{Department of Mathematics and Statistics, University of Jyv\"askyl\"a}
\author[Mazzeo]{Rafe Mazzeo}
\address{Department of Mathematics, Stanford University}

\newtheorem{theorem}{Theorem}
\newtheorem{proposition}[theorem]{Proposition}
\newtheorem{lemma}[theorem]{Lemma}
\newtheorem{corollary}[theorem]{Corollary}

\theoremstyle{definition}

\newtheorem{remark}[theorem]{Remark}

\DeclareMathOperator{\Ker}{Ker}

\newcommand*{\R}{\mathbb{R}}
\newcommand{\RR}{\mathbb R}

\newcommand{\der}{\mathrm{d}}
\newcommand{\eps}{\varepsilon}
\newcommand{\bdf}{x}

\newcommand{\Order}{\mathcal{O}}

\newcommand{\Vol}{\mathrm{Vol}}
\newcommand{\calV}{\mathcal V}
\newcommand{\calU}{\mathcal U}
\newcommand{\calH}{\mathcal H}

\newcommand{\gradh}{\nabla^{\calH}}
\newcommand{\gradv}{\nabla^{\calV}}

\newcommand{\abs}[1]{\left\vert#1\right\vert}
\newcommand{\nabs}[1]{\vert#1\vert}

\newcommand{\norm}[1]{\left\lVert#1\right\rVert}
\newcommand{\nnorm}[1]{\lVert#1\rVert}

\newcommand{\ip}[2]{\left\langle#1,#2\right\rangle}
\newcommand{\nip}[2]{\langle#1,#2\rangle}

\newcommand{\niip}[2]{(#1,#2)}

\newcommand{\iM}{M^\circ}

\newcommand{\del}{\partial}
\newcommand{\gbar}{\overline{g}}
\newcommand{\calC}{\mathcal C}
\newcommand{\calO}{\mathcal O}

\newcommand{\dist}{\mathrm{dist}\,}

\begin{document}

\maketitle

\begin{abstract}
On gas giant planets the speed of sound is isotropic and goes to zero at the surface. Geometrically, this corresponds to a Riemannian manifold whose metric tensor has a conformal blow-up near the boundary. The blow-up is tamer than in asymptotically hyperbolic geometry: the boundary is at a finite distance.

We study the differential geometry of such manifolds, especially the asymptotic behavior of geodesics near the boundary. We relate the geometry to the propagation of singularities of a hydrodynamic PDE and we give the basic properties of the Laplace--Beltrami operator. We solve two inverse problems, showing that the interior structure of a gas giant is uniquely determined by different types of boundary data.
\end{abstract}


\section{Introduction}

The study of propagation of acoustic waves on a gas giant planet leads to a Riemannian geometry that lies between asymptotically hyperbolic geometry and standard geometry with boundary.
Some of the phenomena in this geometry are unlike those seen at either end; for example, constant curvature is not possible.
We set out to study this geometry, the related analytic model, and inverse problems for determining the geometry from boundary measurements.

On a gas giant planet, unlike a rocky planet, the speed of sound goes to zero at the boundary.
Geometrically, the power law decay rate of the speed of sound corresponds to a specific conformal blow-up rate of a Riemannian metric.
This rate is slower than on asymptotically hyperbolic manifolds and the boundary is at a finite distance from interior points.

We study the basic geometry of gas giant Riemannian metrics, including
properties of geodesics near the boundary (Propositions~\ref{prop:F-diffeo} and~\ref{prop:locally-unique-geodesic}),
the Hausdorff dimension of the boundary (Proposition~\ref{prop:bnd-hausdorff-dimension}), and
discreteness of the spectrum of the Laplace--Beltrami operator (Proposition~\ref{prop:Fredholm-and-discrete-spectrum}).

We solve two inverse problems for simple gas giants planets, proving that the metric is uniquely determined by its boundary distance data
(Theorem~\ref{thm:bdy-distane-function}) and that the geodesic X-ray transform is injective (Theorem~\ref{thm:injectivity-in-simple-geometry}).

A brief introduction to gas giant physics and how it leads to our geometric model is given in Section~\ref{sec:physics-intro} below and a more 
detailed model is discussed in Section~\ref{sec:gas-giant-physics}.

\subsection{Gas giant geometry}

Let $M$ be an $(n+1)$-dimensional compact manifold with boundary. A metric $g$ on $M$ is called a 
gas giant metric of order $\alpha \in (0,2)$ if it can be written in the form
\begin{equation}
g = \frac{\overline{g}}{x^\alpha}
\label{gasgiantmetric}
\end{equation}
where $\overline{g}$ is a smooth non-degenerate metric on $M$, including up to its boundary.
Observe that any such metric is incomplete. There are two limiting cases: when $\alpha \to 0$,
$g$ becomes the ordinary incomplete metric $\overline{g}$ on $M$, but when $\alpha \to 2$,
then $g$ converges to a complete asymptotically metric of a type often called conformally compact; cf. e.g.~\cite{Maz1988}. 
We shall typically use a useful normal form.  We may choose local coordinates $(x,y)$ 
on $M$, where $x \geq 0$ and $x = 0$ but $dx \neq 0$ on $\partial M$ and $y$ restricts to 
a coordinate system on the boundary.  There is an associated collar neighborhood of the boundary $\mathcal U
\cong [0,1)_x \times \partial M$ and a smooth family of metrics $h_x$ on $\partial M$ such that
\begin{equation}
g = \frac{ dx^2 + h_x}{x^\alpha}.
\label{gasgiantmetric2}
\end{equation}
This is an analogue of the Graham-Lee normal form for conformally compact metrics. We establish this below
in Section 2.1.

Our goals in this paper are to develop a number of facts about the geometry and analysis of this class of
singular metrics.   The first steps involve a series of calculations concerning the more elementary 
geometric considerations.  We also consider the somewhat more subtle problem of understanding
the asymptotics of escaping geodesics, and of the limiting dynamics of the geodesic flow.
This leads to a first sort of inverse question: is there a way to characterize a gas giant metric intrinsically?
More specifically, if $(M^\circ,g)$ is an open manifold with an incomplete metric, then is it possible
to determine from this metric alone the compactification $M$, as a smooth manifold with boundary,
the metric $\overline{g}$, the constant $\alpha$ and the boundary defining function $x$?  

We consider some deeper inverse problems related to this class of metrics.
In particular, we prove that the
X-ray transform $I_g$ on $(M,g)$ is injective.
In the final sections of this paper we also consider the
Laplace--Beltrami operator $\Delta_g$.
We study its spectrum, mapping properties and whether it is
essentially self-adjoint.

This paper is an initial foray into the analysis and geometry of gas giant metrics. Our aim here is to develop
a number of fundamental results, either ab initio or as consequences of other related studies, which will
then make it possible to consider some deeper inverse problems for this class of metrics.   This paper splits
into two not altogether distinct sections. In the first we develop a number of fundamental facts about the
Riemannian geometry, including the behavior of  geodesics, for gas giant metrics.
Some properties are slightly simpler in the special case $\alpha=1$ but we present all our results for all values $\alpha\in(0,2)$.
The second part of the 
paper studies various analytic properties of the scalar Laplace--Beltrami operators for such metrics. 
In between these two parts, we also prove some Pestov-type identities, which involve the vector field
generating geodesic flow on the cosphere bundle, and use these to solve an inverse problem.

\subsection{Geometry from the equation of state}
\label{sec:physics-intro}

As a leading order approximation, we take a gas giant planet to be a ball and assume all physical quantities to be invariant under rotations.
Spherical symmetry is irrelevant for the geometric model introduced above, but it makes physics simpler.


Many celestial bodies are modelled to leading order as polytropes, a far more detailed discussion of which can be found in~\cite{Hor2004}.
The defining feature of a polytrope is the polytropic equation of state
\begin{equation}
p
=
K\rho^{1+1/n}
\end{equation}
relating the pressure $p$ and the density $\rho$ via the polytropic constant $K$ and the polytropic index $n$.
The leading order approximation to a self-gravitating and spherically symmetric polytropic body can be written in terms of an auxiliary radial function $\theta(r)$ that satisfies
$p(r)=p_0\theta(r)^{n+1}$
and
$\rho(r)=\rho_0\theta(r)^{n}$.
If the ambient dimension is $d$ and the polytropic index satisfies $n>-1$, the function $\theta$ satisfies the Lane--Emden equation
\begin{equation}
\theta''(r)
+
(N-1)r^{-1}\theta'(r)
+ Cr^n
=0,
\end{equation}
where $C>0$.
By rescaling the radial variable one can achieve $C=1$.

At the surface of the body where $r=R$ we have $\theta(R)=0$, and by virtue of being a positive (inside the body) solution to the second order Lane--Emden equation the function $\theta$ must satisfy $\theta'(R)<0$.

The speed of sound can be computed as the (isentropic) derivative
\begin{equation}
c
=
\sqrt{\frac{\partial p}{\partial\rho}}
=
K'\rho^{1/2n}
=
K''p^{1/2(n+1)}
=
K'''\theta^{1/2}
\end{equation}
for new constants $K'$ and $K''$ and $K'''$.
This means that the speed of sound is comparable to the square root of the distance to the surface, no matter the value of the polytropic index.
For gas giants the polytropic index is usually taken to agree with the adiabatic index, which is $n=5/3$ in the case of a monoatomic gas.

The polytropic model is only a leading order approximation and is not expected to hold perfectly.
Bodies are also not perfectly rotationally symmetric due to rotation and inhomogeneities.
Therefore we do not take the polytropic model as the truth, but as a guide to choosing a realistic mathematical model.

If $e$ is the Euclidean Riemannian metric on a smooth domain $B\subset\R^n$, then the speed of sound $c(r)$ can be modeled by the conformally Euclidean Riemannian metric $g=c^{-2}e$.
For a symmetric planet $B$ would be a ball.
If $\bdf$ is a boundary defining function for $B$ (i.e. $\bdf(z)>0$ for $z\in B$, $\bdf(z)=0$ for $z\in\partial B$, and $\der\bdf\neq0$ at $\partial B$), the polytropic model suggests that $c(z)\approx \bdf(z)^{1/2}$, and this is the simplest model for a gas giant.
For a rocky planet the speed of sound has a non-zero limit at the boundary and so $c(z)\approx 1$.

Therefore we take for a general model a speed of sound $c(z)\approx \bdf(z)^{\alpha/2}$.
For a gas giant we expect the value of the parameter $\alpha$ to be $1$ and for rocky planets $0$.
For realistic gaseous celestial bodies we may thus reasonably expect that $\alpha$ is close to $1$.
We thus allow $\alpha\in(0,2)$.
The extreme case $\alpha=0$ corresponds physically to solid bodies and mathematically to manifolds with boundary, and the other extreme $\alpha=2$ corresponds to asymptotically hyperbolic geometry but is far from all planetary models.

Therefore we say that a gas giant metric on a smooth manifold $M$ with boundary is a Riemannian metric $g$ on $\iM$ so that $g=\bdf^{-\alpha}h$, where $\bdf$ is a boundary defining function for $M$ and $h$ is a well-defined Riemannian metric up to the boundary.
The fact that $h$ is neither zero nor infinite at $\partial M$ implies a specific blow-up rate for $g$ near the boundary.
This conformal power-law blow-up is the key geometric feature of gas giant metrics.
Both extremes $\alpha=0$ and $\alpha=2$ are quite well understood mathematically, but the intermediate cases $\alpha\in(0,2)$ have been studied far less.
The physically most relevant case $\alpha=1$ does not appear to be geometrically substantially different from other values in the range we allow apart from some minor conveniences and inconveniences that are not important for the present paper.

For a more detailed physical model for the hydrodynamics of a gas giant planet, see Section~\ref{sec:gas-giant-physics} below.
\subsection*{Acknowledgements}
MVdH was supported by the Simons Foundation under the MATH + X program,
the National Science Foundation under grant DMS-2108175, and the corporate members of the Geo-Mathematical Imaging Group at Rice University. A significant part of the work of MVdH was carried out while he was an invited professor at Centre Sciences des Donn\'{e}es at Ecole Normale Sup\'{e}rieure, Paris. JI and AK were supported by the Research Council of Finland (Flagship of Advanced Mathematics for Sensing Imaging and Modelling grant 359208; Centre of Excellence of Inverse Modelling and Imaging grant 353092; and other grants 351665, 351656, 358047).
AK was supported by the Finnish Academy of Science and
Letters.

\section{The geometry of gas giant manifolds}

We begin with an `extrinsic' study of the metric $g = x^{-\alpha} \gbar$.   Namely, we assume that the metric
takes this form and proceed to study its various geometric properties.  

\subsection{Normal forms and asymptotic curvatures}

A first observation is that if a metric $g$ is know to be a gas giant metric for some $\alpha$, then
this value can be determined from the intrinsic geometry of $g$.

\begin{proposition}
Suppose that $g$ is an $\alpha$-gas giant metric on the interior of some manifold with boundary $M$.
Then $g$ is incomplete, and there is a smoothly varying orthonormal basis of sections for $TM$ such
that the sectional curvatures for $2$-planes spanned by pairs of these basis vectors are asymptotic to
\[
-\frac{2\alpha}{(2-\alpha)^2}\dist(\cdot, \del M)^{-2},   \ \ \  \mbox{and} \ \ \ 
-\frac{\alpha^2}{(2-\alpha)^2} \dist(\cdot, \del M)^{-2}.
\]
\label{curvasym}
\end{proposition}

Thus $\alpha$ can be recovered from these asymptotic sectional curvatures.  

We prove this Proposition below, but before doing so, first describe a ``normal form'' for the metric near 
the boundary. This is modelled on a very useful normal form, due to Graham and Lee~\cite[Lemma 5.2]{GL1991}, in the case when 
$\alpha = 2$, in which case the metric $g$ is complete, and is called conformally compact. In that 
case, one can define $\gbar = x^2 g$ where $x$ is any choice of boundary defining function, and by
definition, $\gbar$ is a smooth non-degenerate metric up to the boundary. The restriction of $\gbar$ 
to $\del M$ is a metric on the boundary; however, replacing $x$ by $x' = ax$ where $a$ is any positive
smooth function results in a new metric on $\del M$ conformal to the first one. In other words, only
the conformal class of the metric is well-defined.  The Graham--Lee theorem states that if $h_0$ is
any representative of that conformal class, there is a unique boundary defining function $x$ such that
\[
g = \frac{dx^2 + h(x, y, dy)}{x^2}, \quad h(0,y, dy) = h_0.
\]
Here $h$ is a family of metrics on $\del M$ (pulled back to the level sets $x = \mathrm{const.}$) 
depending smoothly on $x$, and $y$ is any local coordinate system of the boundary.  In particular, 
$-\log x$ is a distance function for the metric $g$.

In the gas-giant setting we can attempt to prove the same thing, but there is no longer  ``free data'' (analogous 
to the choice of representative of the conformal class).
\begin{proposition}
\label{prop:adapted-coordinates}
Let $g$ be an $\alpha$ gas giant metric. Then there is a well-defined metric $h_0$ on $\del M$, and an associated 
boundary defining function $x$ on $M$ such that
\[
g = \frac{dx^2 + h(x, y, dy)}{x^\alpha}, \quad \mbox{where}\qquad h(0,y,dy) = h_0.
\]
\label{normalform}
\end{proposition}
\begin{proof}  First choose an arbitrary boundary defining function $\tilde{x}$. We modify it in two steps. In
the first, we seek a new boundary defining function $\hat{x} = a \tilde{x}$ such that $|d\hat{x}/\hat{x}^{\alpha/2}|^2_g|_{\del M} \equiv 1$.
For this, we compute
\[
\frac{ d(a \tilde{x} )}{ (a \tilde{x})^{\alpha/2}} =  a^{1-\alpha/2} \frac{d\tilde{x}}{ \tilde{x}^\alpha/2} + \calO(\tilde{x}),
\]
hence we simply need choose $a$ along $\del M$ so that $ a^{2-\alpha} | d\tilde{x}|^2/ \tilde{x}^\alpha \equiv 1$ there.

The metric $h_0$ on $\del M$ is then defined as the pullback of $\hat{x}^\alpha g$ to the boundary. The computation
above shows that there is no leeway: the $1$-jet of the boundary defining function, and hence this boundary metric,
are completely fixed by the requirement that $|dx/x^\alpha|_g \equiv 1$. 

We now make a further change, setting $x = e^{\omega} \hat{x}$, and study the equation $| dx/x^{\alpha/2}|^2_g \equiv 1$,
not just at the boundary but in the collar neighborhood of the boundary. Writing $\hat{g} = \hat{x}^{\alpha} g$, we can
rewrite this as
\[
e^{(2-\alpha) \omega}  \frac{ |d\hat{x} + \hat{x} d\omega|_g^2}{ \hat{x}^\alpha} = e^{(2-\alpha)\omega} |d\hat{x} + \hat{x} d\omega|_{\hat{g}}^2 = 1.
\]
Expanding and rearranging yields
\[
|d\hat{x}|^2_{\hat{g}} + 2 \hat{x} \langle d\hat{x} , d \omega\rangle_{\hat{g}} + \hat{x}^2 |d \omega|^2_{\hat g} = e^{(\alpha-2)\omega}. 
\]
Using the normalization of $\hat{x}$ and writing $\langle d\hat{x}, d\omega \rangle_{\hat{g}} = \del_{\hat{x}} \omega$, we recast this in the form
\begin{equation}
\hat{x} \del_{\hat{x}} \omega = - \hat{x}^2 |d\omega|^2 + (1-|d\hat{x}|^2)  +  G(\omega) \omega,
\label{fhj}
\end{equation}
where $G(\omega) = \omega^{-1}( e^{(\alpha-2)\omega} - 1)$ is a smooth function of $\omega$ (including where $\omega$ vanishes). It is important
to note that $G(0) = \alpha-2 < 0$. 

This is a characteristic Hamilton--Jacobi equation.   Fortunately the main result in~\cite{GK2012} is an existence theorem for equations of precisely this form.
That theorem applies to equations of the form
\[
\hat{x} \del_{\hat{x}} \omega = F(x, y, \omega, \del_y \omega), \quad \omega(0,y) = \omega_0(y),
\]
where $F(x,y,\omega, q)$ is smooth and satisfies 
\[
F(0,y,\omega_0, \del_y \omega_0) = 0,\ F_\omega(0,y,\omega_0, \del_y \omega_0) < 1,\ \ F_q(0,y,\omega_0, \del_y \omega_0)  = 0.
\]
The conclusion in~\cite{GK2012} is that there exists a unique smooth solution in some small interval $0 \leq \hat{x} < \hat{x}_0$.   
In the proof they observe that the stronger condition $F_\omega < 0$ at $(0,y,\omega_0, \del_y \omega_0)$ implies that the solution is unique even
amongst continuing solutions.

To apply this theorem to our setting, we impose the initial condition $\omega(0,y) = \omega_0 = 0$.  We then write
$\hat{x}\del_{\hat{x}}\omega + \hat{x}^2 |d\omega|^2_{\hat{g}} = H(\hat{x}\del_{\hat{x}}\omega, \hat{x}\del_{y}\omega)$, where $H(q_1, q_2)$
satisfies $H_{q_1}(0,0) = 1$, $H_{q_2}(0,0) = 0$. Applying the implicit function theorem, we can thus rewrite~\eqref{fhj} as
\[
\hat{x} \del_{\hat{x}}\omega = \mathcal F(\hat{x}, y, \omega, \del_y \omega)
\]
where the differential of $\mathcal F$ in its third argument at $\omega_0 = 0$ equals $G(0)$. Since $F_\omega G(0) = \alpha-2 < 0$, the result of
Graham and Kantor can be applied. In fact, even the stronger form, which gives uniqueness even amongst all continuous solutions, also holds. 
\end{proof}

We now return to the assertion about curvature asymptotics.

\begin{proof}[Proof of Proposition \ref{curvasym}]
We first compute sectional curvatures for the warped product metric $ x^{-\alpha}( dx^2 + h_0)$.
The shortest way uses Cartan's method of moving frames, which we recall briefly. 
We choose a $g$-orthonormal family of $1$-forms $\{\omega_i\}$ which span $T_p^* M$ at each point.
Thus, essentially by definition, $g = \sum \omega_i \otimes \omega_i$.  A simple lemma states 
that there exist uniquely defined $1$-forms $\omega_{ij}$ which are skew-symmetric in the indices,
i.e., $\omega_{ji} = -\omega_{ij}$, such that
\[
d\omega_i = \sum_j \omega_{ij} \wedge \omega_j.
\]
This is called Cartan's lemma, and the forms $\omega_{ij}$ encode the Levi-Civita connection. We then define $2$-forms
\[
\Omega_{ij} \coloneqq d\omega_{ij} - \sum_k \omega_{ik} \wedge \omega_{kj}.
\]
It is then not difficult to show (and this is explained in many sources) that
\[
\Omega_{ij} = - \sum_{k, \ell} R_{ijk\ell} \omega_k \wedge \omega_\ell,
\]
where $R_{ijk\ell}$ are the components of the Riemann curvature tensor in this basis at each point. 

We apply this as follows. Let $\bar{\omega}_\beta$, $\beta = 1, \ldots, n-1$, denote a smoothly varying orthonormal basis of
$1$-forms on $(\del M, h_0)$, and write
\[
\omega_0 = \frac{dx}{x^{\alpha/2}}, \quad \omega_\beta = \frac{\bar{\omega}_\beta}{x^{\alpha/2}}.
\]
A short calculation then shows that for $\beta, \gamma = 1, \ldots, n-1$, 
\[
\omega_{\beta \gamma} = \bar{\omega}_{\beta \gamma}, \qquad \omega_{\beta 0} = \frac{\alpha}{2} x^{(\alpha-2)/2} \omega_\beta.
\]
Here $\bar{\omega}_{\beta \gamma}$ are the connection $1$-forms for the metric $h_0$ on $\del M$ (extended to the
neighborhood $\calU$ by the product decomposition). 

Finally we compute that
\[
\Omega_{\beta \gamma} = \frac{\alpha^2}{4} x^{\alpha-2} \omega_\beta \wedge \omega_\gamma + \calO(x^{\alpha}), 
\quad \mbox{and}\quad
\Omega_{\beta 0} = \frac{\alpha}{2} x^{\alpha-2} \omega_\beta \wedge \omega_0 + \calO(x^{\alpha-1}).
\]
The estimate of the remainder term uses, for example, that $|\Omega_{\beta \gamma}|_g^2 = x^{2\alpha} |\Omega_{\beta \gamma}|^2_{\gbar}$.
We conclude that the principal components of the curvature tensor (which agree with the corresponding
sectional curvatures because of our use of orthonormal coframes) satisfy
\[
R_{\beta 0 \beta 0} \sim -\frac{\alpha}{2} x^{\alpha-2}, \quad R_{\beta \gamma \beta \gamma} \sim - \frac{\alpha^2}{4} x^{\alpha-2}.
\]

The function $x$ is related to the distance function $s$ by $(1-\alpha/2) x^{1-\alpha/2} = s$, hence
\[
R_{\beta 0 \beta 0} \sim - \frac{ -2\alpha}{ (2-\alpha)^2}, \quad R_{\beta \gamma \beta \gamma} \sim -\frac{\alpha^2}{(2-\alpha)^2},
\]
as claimed. 

We have shown that any gas giant metric can be written in this simple warped product form up to remainders which are $\calO(x)$.
However, there is something mildly circular in that we used an initial knowledge of $\alpha$ in proving that normal form.
To show that this is not a true issue, observe that we can carry out with only moderately more work the same computations
as above if we only know that the metric $g$ is a gas-giant metric for some parameter $\alpha$, and have set $\gbar = x^{\alpha} g$
for an arbitrary boundary defining function $x$.  The leading asymptotics then determine the value of $\alpha$ just as above.
\end{proof}

We list a few more basic properties.

\begin{proposition}
If $(M,g)$ is a gas giant metric, then $\mathrm{Vol}\,(M,g) < \infty$ if and only if $\alpha < 2/n$. If $\alpha > 2/n$,
then $\mathrm{Vol}\,(\{x \geq \eps\}) \sim C \eps^{1-n\alpha/2}$, while if $\alpha = 2/n$, then 
$\mathrm{Vol}\,(\{x \geq \eps\}) \sim -C \log \eps$.
\label{volgrowth}
\end{proposition}
\begin{proof}
In the special coordinates above, $dV_g = x^{-n\alpha/2} dx dV_h$, so the total volume is finite if $-n\alpha/2 > -1$, i.e., 
$\alpha < 2/n$. The other assertions are immediate. 
\end{proof}

\begin{proposition}
The second fundamental form of the level sets $\{x = \eps\}$ are strictly convex.
\label{2ndfundform}
\end{proposition}
\begin{proof}
This is a standard computation, which is left to the reader. The conclusion is that 
\[
\nabla_{\del_{y_i}} (-x \del_x) = \frac{\alpha}{2} \del_{y_i} + \calO(x).
\]
This shows that the second fundamental form of these level sets is, asymptotically, $\alpha/2$ times the identity, and in particular
is positive definite. 
\end{proof}

\subsection{Geodesics}
We now turn to a study of the geodesic flow on $(M,g)$. 

In the following we always use an adapted coordinate system $(x,y)$ near the boundary, where $x$ is the special boundary defining
function obtained in Proposition~\ref{normalform} and $y$ is any coordinate system on the boundary.  We denote by $(\xi, \eta)$ 
the associated covectors.  We shall use the Hamiltonian formalism, namely we write the equations for the bicharacteristics 
for the Hamiltonian function 
\[ 
H(x,y,\xi,\eta) = \frac12 | (\xi,\eta)|^2_{g_{x,y}} = \frac12 (x^{\alpha} \xi^2 + x^\alpha h^{ij}(x,y) \eta_i \eta_j).
\]
These bicharacteristics are curves in $T^* M$ which project to the geodesics on $M$.    These equations are:
\begin{equation}
\label{eqn:coflow}
\begin{split} 
\dot{x} & = \frac{\del H}{\del \xi} = x^\alpha \xi, \quad \dot{y_i} = \frac{\del H}{\del \eta_i} = x^\alpha h^{ij}(x,y) \eta_j  \\
\dot{\xi} &= -\frac{\del H}{\del x} =  - \alpha x^{-1} H(x,y,\xi,\eta) - \frac12 x^\alpha \frac{\del h^{ij}}{\del x} \eta_i \eta_j, \\
\dot{\eta_i} & = -\frac{\del H}{\del y_i} = -\frac12 x^\alpha \frac{\del h^{jk}(x,y)}{\del y_i} \eta_j \eta_k.
\end{split}
\end{equation}
We may as well restrict to geodesics of a fixed speed, and thus suppose that $H \equiv 1/2$ along the solution curves. This simplifies the first
summand in the equation for $\dot{\xi}$ to being simply $-\alpha/2x$.   We often write a bicharacteristic as $(z(t), \zeta(t))$, where
$z(t) = (x(t), y(t))$ and $\zeta(t) = (\xi(t), \eta(t))$. 

Before we begin to analyze this system, there are some preliminary observations.  First, $ x^\alpha (\xi^2 + h^{ij} \eta_i \eta_j) \equiv 2$ along
each orbit, so from this it follows that if $A^{ij}$ is any matrix which is uniformly bounded on $M$, e.g., one written in terms of partial derivatives of the $h^{ij}$ 
with respect to any of the variables $x$ or $y_k$, then 
\begin{equation}
\left| x^\alpha A^{ij} \eta_i \eta_j \right| \leq C,
\label{bd1}
\end{equation}
along each orbit, where $C$ depends only on the norm of $A$. In the following, we use $\calO(1)$, $\calO(x^\alpha)$, etc., to denote quantities 
which are bounded by $C$, $C x^\alpha$, etc., where the constants $C$ depend only on the metric and are independent of the orbit. 

\begin{lemma}
\label{lma:monotonous-towards-b}
For $\eps > 0$ small enough, if $\gamma(t) = (z(t), \zeta(t))$ is any bicharacteristic with $x(0) < \eps$ and
$\xi(0) = 0$, then $\xi(t) < 0$ for all $t \in \RR$. 
\label{lemma1}
\end{lemma}
\begin{proof}
The hypothesis is invariant with respect to replacing $t$ by $-t$, so we prove the assertion for $t \geq 0$.  First observe that, by~\eqref{bd1},
\[
\dot{\xi} = -\alpha x^{-1} + \calO(1) < -\frac12 \alpha \eps^{-1}  < 0
\]
if $\eps$ is sufficiently small.  Again, the penultimate inequality here is independent of the trajectory.  

This argument shows that if $\xi(0) = 0$, then $\xi(t) < 0$ for $t > 0$ sufficiently small, but in fact it shows that for any $t_0 > 0$, if $\xi(t_0) < 0$ and $x(t_0) <
\eps$, then $\xi(t)$ remains bounded above by a strictly negative constant.  This proves that $\xi(t) < 0$ for all $t \geq 0$, and for any $t_0 > 0$,
$\xi(t) \leq -c < 0$ for $t \geq t_0$. 
\end{proof}

\begin{lemma}
\label{lma:finite-exit-time}
If $\gamma(t) = (z(t), \zeta(t))$ is any bicharacteristic with $x(0) < \eps$, where $\eps$ is chosen as in Lemma~\ref{lemma1}, 
and $\xi(0) \leq 0$, then $z(t)$ converges to a unique point $(0, \bar{y}) \in \del M$ at some finite time $T > 0$ and $\eta(t)$ converges
to some $\bar{\eta}$ as $t \nearrow T$ as well. 
\label{lemma2}
\end{lemma}
\begin{proof}
We have just shown that the function $x(t)$ is strictly monotone decreasing.  Denote the maximal time of existence by $T \leq \infty$. 
There are a number of possibilities: either $T < \infty$ or $T = \infty$, and in each of these cases, either $x(t) \searrow x_0 > 0$ as $t \nearrow T$
or else $x(t) \searrow 0$. We aim to show that $T < \infty$ and $x(t) \searrow 0$.   

Suppose first that $x(t) \searrow x_0 > 0$. If in addition $T < \infty$, then the system of equations remains non-degenerate and we could simply 
take a limit as $t \to T$ to define $\gamma(T)$ and then continue the solution for later times $t > T$.  On the other hand, if $T = \infty$, then
using that $\dot{\xi}(t) \leq -c < 0$ for $t \geq t_0$, we obtain $\xi(t) \to -\infty$, which would contradict that $x_0^\alpha \xi^2 < x^\alpha \xi^2 \leq 1$.   
Neither of these scenarios are possible, hence $x(t) \searrow 0$.

We next show that $\gamma(t)$ reaches $x =0$ in finite time. Since $x$ is monotone, we may use it as the independent parameter. Thus, writing
$\xi = \xi(x)$, we have
\[
\frac{d\xi}{dx} = \frac{ - (\alpha/2) x^{-1} + \calO(1)}{ x^\alpha \xi} \Longrightarrow  \frac{d\,}{dx} \xi(x)^2 = - \alpha x^{-\alpha-1} + \calO(x^{-\alpha}).
\]
Writing the final term as $x^{-\alpha} F$, where $F$ is bounded, and integrating from $x$ to $1$, gives
\[
\xi(1)^2 - \xi(x)^2 = (1 - x^{-\alpha}) + \int_x^1 s^{-\alpha} F(s)\, ds,
\]
whence $\xi(x) = - x^{-\alpha/2}(1 + \calO(x) + \calO(x^\alpha) )$. (The case $\alpha = 1$ is of course slightly different, but we omit the details.) 
Now insert this into the equation for $\dot{x}$ to get that 
\begin{equation}
\begin{split}
\frac{dx}{dt} = x^\alpha \xi = & - x^{\alpha/2}(1 + \calO(x) + \calO(x^\alpha) ) \Longrightarrow \\ 
& x^{-\alpha/2} \dot{x} = - 1 + \calO(x) + \calO(x^\alpha). 
\end{split}
\label{morerefeqn}
\end{equation}
Bounding the last two terms by $C (\eps + \eps^\alpha)$, and integrating from $t_0$ to $t_1$, we get
\[
x(t_1)^{1-\alpha/2} = x(t_0)^{1-\alpha/2}  - (1-\alpha/2)(t_1 - t_0) + \calO(\eps + \eps^\alpha) (t_1 - t_0).
\]
As $t_1 \nearrow T$, $x(t_1)^{1-\alpha/2} \to 0$, which then shows that $t_1$ cannot become arbitrarily large. This proves that $T < \infty$.

We next observe that by the Hamiltonian constraint, $\dot{\eta}_i = \calO(1)$, and and thus $\eta_i(t)$ converges to some limiting value $\bar{\eta}_i$ as 
$t \to T$ since $T$ is finite.  Using this, we also conclude that $y_i(t) \to \bar{y}_i$, and furthermore that
\begin{equation}
(y(t), \eta(t)) = (\bar{y}, \bar{\eta}) + \calO(x^\alpha).
\label{firstestyeta}
\end{equation}

Note, however, that $\xi(t)$ is unbounded, and more specifically,  $\xi(t) \sim - x(t)^{-\alpha/2} \to -\infty$ as $t \nearrow T$.
\end{proof}

We now improve these estimates by showing that along a fixed trajectory, the functions $x(t)$, $y(t)$, $\xi(t)$ and $\eta(t)$ have
complete asymptotic expansions in powers of $\tau = T - t$ as $\tau \to 0$. This is achieved by an iteration argument and a careful 
examination of the methods used in the preceding proof. To simplify notation below, we use $\tau$ as a new independent
variable, and for any function $f(\tau)$, denote $df/d\tau$ by $f'$ (so $f' = - \dot{f}$).   We proceed with the calculations,
and summarize the outcomes of all of this at the end. 

First, integrate $x^{-\alpha/2} x' = \calO(1)$ from $0$ to $\tau$ to get $x(\tau) = \calO( \tau^{\frac{2}{2-\alpha}})$. Substituting this into
\eqref{morerefeqn} yields $x^{-\alpha/2} x' = 1 + \calO( \tau^{2/(2-\alpha)} + \tau^{2\alpha/(2-\alpha)})$, which then implies that
\begin{equation}
x(\tau) = (1 - \alpha/2)^{2/(2-\alpha)} \tau^{2/(2-\alpha)} (1 + \calO(\tau^{2/(2-\alpha)} + \tau^{2\alpha/(2-\alpha)}) ).
\label{expx}
\end{equation}
This gives a leading asymptotic term for the function $x(\tau)$. 

For the next step, observe that since we have already proved that $y$ and $\eta$ remain bounded, the equations of motion show that
$(y', \eta')  = \calO(x^\alpha) = \calO( \tau^{2\alpha/(2-\alpha)})$, so that
\begin{equation}
(y(\tau), \eta(\tau)) = (\bar{y}, \bar{\eta}) + \calO( \tau^{ (2 + \alpha)/(2-\alpha)}).
\label{expyeta}
\end{equation}
Finally, 
\begin{equation}
\xi(\tau) = - (1-\alpha/2)^{2/(2-\alpha)} \tau^{-\alpha/(2-\alpha)} (1 + \calO( \tau^{2/(2-\alpha)} + \tau^{2\alpha/(2-\alpha)}) ). 
\label{expxi}
\end{equation}

The equations~\eqref{expx}, \eqref{expyeta} and~\eqref{expxi} show that each of the functions $x(\tau), y(\tau), \xi(\tau), \eta(\tau)$ has
a leading asymptotic term plus a lower order remainder as $\tau \to 0$.  For many purposes this is sufficient. However,  it is straightforward
to set up an inductive scheme to prove the existence of complete polyhomogeneous expansions for these functions in powers of $\tau$. 
(When $\alpha = 1$, these expansions also involve positive integer powers of $\log \tau$ as well.) This is done by iteratively substituting
the partial expansions of these functions into the equations of motion and integrating from $0$ to $\tau$, which produces an expansion 
with one further term in the asymptotic plus an error term which vanishes even more quickly. 

Since it will be very helpful below, we carry out the first step of this iteration. In the following, set
\[
c_\alpha = (1-\alpha/2)^{2/(2-\alpha)},
\]
and for simplicity, indicate higher order remainders by ``$...$''.  Now, insert the expansions~\eqref{expx} and~\eqref{expyeta} into the 
equation for $y_i'(\tau)$ to get that
\begin{multline*}
y_i'(\tau) =  x^\alpha h^{ij}(x,y) \eta_j \\
= ( c_\alpha \tau^{2/(2-\alpha)} + \ldots)^\alpha h^{ij}( c_\alpha \tau^{2/(2-\alpha)} + \ldots,
\bar{y} + \ldots) (\bar{\eta}_j + \ldots) \\
= c_\alpha^\alpha \tau^{2\alpha/(2-\alpha)} h^{ij}(0,\bar{y}) \bar{\eta}_j + \ldots,
\end{multline*}
Note that $h^{ij}(0,\bar{y}) \bar{\eta}_j = \bar{v}^i$ is the $i^{\mathrm{th}}$ coordinate of the vector $\bar{v}$ which is $h_0$-dual to $\bar{\eta}$
at $\bar{y}$.  Thus
\begin{equation}
y(\tau) = \bar{y} + c_\alpha' \tau^{(2+\alpha)/(2-\alpha)} \bar{v} + \ldots,
\label{y2}
\end{equation}
where $c_\alpha' = \alpha^{-1}( (2-\alpha)/2)^{(2+\alpha)/(2-\alpha)}$. 

From this we immediately deduce the following.

\begin{corollary}
Any geodesic $(x(t), y(t))$ which approaches the boundary does so along a curve asymptotic to
\[
y - \bar{y} =  c_\alpha'' x^{(2+\alpha)/2} \bar{v}
\]
for some $\bar{v} \in T_{\bar{y}} \del M$, where $c_\alpha''$ is a universal constant depending only on~$\alpha$. 
\label{geodas}
\end{corollary}

Collecting all of the calculations, and proceeding as explained above, we have proved the following result.

\begin{proposition}
Each trajectory $(z(t), \zeta(t))$ which remains in the region $\{x < \eps\}$ for $t \geq 0$ reaches the boundary
at $x = 0$ at some finite time $T$.  The coordinate functions $(x(t), y(t), \xi(t), \eta(t))$ for a given trajectory admit complete asymptotic 
expansions in powers of $T-t$ (and when $\alpha=1$, also $\log (T-\tau)$. In particular, $(y(t), \eta(t))$ converges to some fixed 
point $(\bar{y}, \bar{\eta})$ in the cotangent bundle of the boundary as $t \to T$.
\label{exptraj}
\end{proposition}

We now consider all points $(z_0, \zeta_0)$, $\zeta_0 \in T_{z_0}^* M$, with $0 < x_0 < \eps$, where the forward trajectory $(z(t), \zeta(t))$ 
remains in the region $x < \eps$ for all $t \geq 0$ and converges to $x=0$. To simplify matters, assume that $\xi_0 = 0$, so that 
$\eta_0$ satisfies the Hamiltonian constraint $x_0^\alpha |\eta_0|^2_{h(x_0,y_0)} = 1$.  Our goal is to determine the dependence of the exit time $T$ and 
exit point $(\bar{y},\bar{\eta})$ as functions of $(x_0,y_0,\eta_0)$. 
\begin{lemma}
\label{lma:exit-time}
The function $T(x_0, y_0, \eta_0)$ is smooth when $x_0 > 0$ and has a complete asymptotic expansion in powers of $x_0$ as $x_0 \to 0$,
with leading term $T \sim c x_0^{1-\alpha/2}$ for some $c > 0$. 
\label{depT}
\end{lemma}
\begin{proof}
Strictly speaking, the analysis in the preceding proof assumes that $\xi(0) < 0$. We arrange this by first using that the one-parameter family
of local diffeomorphisms $\Phi_t$ associated to this flow for some small time $t = \ell(x_0,y_0, \eta_0) > 0$ defines a smooth map 
$\Phi_{t(x_0)} \colon (x_0, y_0, \eta_0) \mapsto (x_1, y_1, \xi_1, \eta_1)$.  We choose this function $\ell(x_0,y_0,\eta_0)$ so that $x_1 = x_0/2$. The ``height'' $x_1$
depends on all the variables $(x_0, y_0, \eta_0)$ (and the function $\ell$ too), so the image of this map as $\eta_0$ varies but $(x_0, y_0)$ remains fixed
is a small $(n-2)$-sphere which is not of constant height, but along which $\xi_1$ is everywhere negative.  By Lemma~\ref{lemma1}, 
the continuing trajectory converges to~$\del M$.  

Now use that $x(\tau) = c \tau^{2/(2-\alpha)} + \ldots$, which implies, equivalently, that $\tau = c' x^{(1-\alpha/2)}  + \ldots$.  These functions (and the
terms in the expansions) are smooth in all the remaining data.  This shows that the time $\tau$  needed to move along this given trajectory from 
$(0, \bar{y})$ to $(x_1, y_1)$ depends smoothly on $(x_0, y_0, \eta_0)$ and is polyhomogeneous in $x_0$.  

We have just shown that elapsed time $h(x_0, y_0, \eta_0)$ for the path to move from $x_0$ to $x_0/2$ is on the order of $c'' x_0^{1-\alpha/2}$ for some $c''>0$.
This function is readily seen to be polyhomogeneous as $x_0 \to 0$, as is the concatenation with the map that gives the elapsed time for the trajectory
to move from $x_0/2$ to the boundary. 
\end{proof}

We next study the ``endpoint mapping'' from the set of initial conditions $\mathcal S \coloneqq \{(x_0, y_0, \eta_0)\colon H(x_0, y_0, 0, \eta_0) = 1/2\}$ to 
the limiting covector on the boundary:
\[
F\colon \mathcal S \longrightarrow T^*\del M, \qquad    F( x_0, y_0, \eta_0) = (\bar{y}, \bar{\eta}).
\]
Of course, $F$ is well-defined only when restricted to the set $\mathcal S_\eps = \{(x_0, y_0, \eta_0) \in \mathcal S: 0 < x_0 < \eps\}$ for some
sufficiently small $\eps$, and we henceforth fix such an $\eps$ and the restriction of $F$ to this set. Note that both the Hamiltonian constraint 
set $\mathcal S_\eps$ and $T^* \del M$ are $(2n-2)$-dimensional. In the following, we systematically identify
covectors $\zeta = (\xi,\eta)$ with vectors $v$ using the metric $g$. However, for covectors $(\bar{y}, \bar{\eta}) \in T^*_{\bar{y}}\del M$,
we identify $\bar{\eta}$ with a vector $\bar{v} \in T_{\bar{y}} \del M$ via the metric~$h_0$.

\begin{proposition}
\label{prop:F-diffeo}
The map $F: \mathcal S_\eps \longrightarrow T^* \del M$ is a diffeomorphism onto its image.  Furthermore, it is smooth,
in a precise sense to be made explicit during the course of the proof, in the limit as $x_0 \to 0$. 
\end{proposition}
\begin{proof}
First note that along geodesics starting on
$\mathcal S_\eps$, we have that $\xi = - x^{-\alpha/2} ( 1/2 - x^\alpha h^{ij}(x,y) \eta_i \eta_j)^{1/2}$. Inserting this into 
the equation for $\dot{x}$ yields that
\[
\frac{dx}{dt} = -x^{\alpha/2} (1/2 - x^\alpha h^{ij}(x,y) \eta_i \eta_j)^{1/2} \eqqcolon - x^{\alpha/2} K.
\]
The quantity $K$ is simply the second factor with the square root. As we have done before, let us shift to using $x$ as the independent variable.  
We can then rewrite the equations for $\dot{y}_i$ and $\dot{\eta}_i$ as
\[
\begin{split}
\frac{dy_i}{dx} & = \frac{ dy_i/dt}{ dx/dt} = - x^{\alpha/2} h^{ij}(x,y) \eta_j,  \\
\frac{d\eta_i}{dx} & = \frac{ d\eta_i/dt}{dx/dt} = \frac12 x^{\alpha/2} \del_{y_i} h^{pq}(x,y) \eta_p \eta_q.
\end{split}
\]

What we have done is to rewrite the equations for $(y,\eta)$ as ``self-contained'' equations involving only the new independent variable
$x$ and $(y,\eta)$. This system takes the form
\[
\frac{d\,}{dx} \begin{bmatrix} y \\ \eta \end{bmatrix} =  x^{\alpha/2} K(x,y,\eta) G(x,y,\eta), \ \  \mbox{where}\ \ 
G(x,y,\eta) = \begin{bmatrix} - h^{ij}\eta_j \\ \frac12 \del_{y_i} h^{pq} \eta_p \eta_q \end{bmatrix}.
\]
Rewrite this as $x^{-\alpha/2}\frac{d\,}{dx}\begin{bmatrix} y \\ \eta \end{bmatrix} = K(x,y,\eta) G(x,y,\eta)$.  This suggests that we reparametrize again,
setting $u = x^{1+\alpha/2}/(1 + \alpha/2)$ so that $\frac{d\,}{du} = \frac{dx}{du} \frac{d\,}{dx} = x^{-\alpha/2} \frac{d\,}{dx}$. The system then becomes
\[
\frac{d\,}{du} \begin{bmatrix} y \\ \eta \end{bmatrix} =  K(x(u), y, \eta), G(x(u), y, \eta).
\]

Finally, the endpoint map we are studying corresponds to the flow of this system between the two values $u_0 = x_0^{1+\alpha/2}/(1+\alpha/2)$
and $u_1 = 0$.   Since $x(u) = (1+\alpha/2)^{2/(2+\alpha} u^{2/(2+\alpha)}$, the functions on the right are polyhomogeneous in $u$, but not
smooth at $u=0$. However, the lack of full regularity in the independent variable is not relevant in the key fact needed here, which is
smooth dependence on initial conditions $(u_0, y_0, \eta_0)$.   This map $\mathcal S_\eps \ni (x_0, y_0, \eta_0) \mapsto  (y(0), \eta(0))$
is this thus smooth, and patently reversible, hence defines a diffeomorphism from the domain $\mathcal S_\eps$ to its image, an open subset
of $T^* \del M$. 

For the final statement, we employ a scaling argument to study this map as $x_0 \to 0$.  Fix a point $(0, 0) \in \del M$, and consider 
the family of dilations $\delta_\lambda: (x,y) \mapsto (\lambda x, \lambda y)$. The pullback of the fixed metric $g$ with respect to $\delta_\lambda$ is
\[
\delta_\lambda^* ( x^{-\alpha} (dx^2 + h_{ij}(x,y) dy^i dy^j) ) = \lambda^{2-\alpha} x^{-\alpha} (dx^2 + h_{ij}(\lambda x, \lambda y) dy^i dy^j),
\]
and after normalizing, this has a limit:
\[
\lim_{\lambda \to 0}  \lambda^{\alpha-2} \delta_\lambda^* g = x^{-\alpha}(dx^2 + h_{ij}(0,0) dy^i dy^j).
\]
This last metric is defined on the entire half-space $\RR^{n}_+ = \{(x,y) \in \RR^n: x > 0\}$.  The (co)geodesic flow of these dilated rescaled
metrics are simply reparametrizations of the geodesics for the initial metric $g$. 

We employ this as follows. To understand the behavior of $F(x_0, y_0, \eta_0)$ as $x_0 \to 0$, it suffices to consider the family of
mappings $F_{x_0}(1, y_0, \eta_0)$ which are defined in the same way as $F$, but for the family of rescaled metric $x_0^{\alpha-2} \delta_{x_0}^* g$. 
These rescaled metrics converge smoothly as $x_0 \to 0$, and this implies easily that this family of mappings $F_{x_0}$ also converge smoothly.
\end{proof}

With this analysis, we can now use the map $F$ to understand further maps of interest.

\begin{proposition}
\label{prop:locally-unique-geodesic}
Let $y_1$ and $y_2$ be two nearby points on $\del M$. Then there exists a unique geodesic $\gamma$ which
connects $y_1$ to $y_2$. 
\end{proposition}

\begin{proof}
Given any $(y_1, \eta_1)$ and a point $(x_0, y_0)$ with $y_0$ sufficiently near to $y_1$, $x_0$ sufficiently small and $|\eta_1| \leq C$, there exists
a unique trajectory $(x(t), y(t), \xi(t), \eta(t))$ with initial condition $(x_0, y_0, 0, \eta_0)$ for some $\eta_0$ satisfying $H(x_0, y_0, 0, \eta_0) = 1/2$
and such that $(y(t), \eta(t)) \to (y_1, \eta_1)$. 
Now follow this trajectory past $(x_0, y_0)$. This continuation hits the boundary at some point $(\bar{y}, \bar{\eta}) = F(x_0, y_0, -\eta_0)$. 
The elapsed time for the entire trajectory is $T(x_0, y_0, \eta_0) + T(x_0, y_0, -\eta_0)$. This defines a smooth invertible map, which is the analogue 
of the scattering relation in this setting, 
\[
E\colon T^* \del M \to T^* \del M, \ \ E(y_1, \eta_1) = (\bar{y}, \bar{\eta}).
\]
Now fix $y_1$ and project $E$ off the $\bar{\eta}$ component. Since there are no conjugate points, this defines a diffeomorphism from a small punctured ball 
$B_c(0)\setminus \{0\} \subset T^*_{y_1} \del M$ of non-zero covectors to its image, a punctured neighborhood of $y_1$ in $\del M$.  Thus to 
every $y_2$ in this neighborhood, there exists some $\eta_1$ such that $E(y_1, \eta_1(y_2)) = (y_2, \eta_2)$ for some $\eta_2$. 
This associates to the pair $(y_1, y_2)$ first the covector $(y_1, \eta_1)$ and then the apex of the corresponding geodesic $F^{-1}(y_1, \eta_1)$, and 
finally the entire geodesic.
\end{proof}

We have now shown that the interior distance function $d_g(y_1, y_2)$ is well-defined for $(y_1, y_2)$ lying in a sufficiently small punctured
neighborhood of the diagonal of $(\del M)^2$.   Let us reparametrize the space of such pairs with the new variables $(\bar{y}, \bar{v})$; here
$\bar{y}$ is defined as the midpoint of the $h_0$-geodesic connecting $y_1$ to $y_2$ and $\bar{v}\in T_{\bar y}\partial M$ is the tangent vector 
to that geodesic (in the direction from $y_1$ toward $y_2$) with length $d_g(y_1, y_2)/2$.  Thus $y_1, y_2 = \exp^{h_0}_{\bar{y}}(\pm \bar{v})$. 
Write $\bar{v} = r \omega$ in spherical coordinates, so $r = d_g(y_1, y_2)/2 \geq 0$ and $\omega \in S^{n-2}$.    
\begin{corollary} 
The interior distance function $d_g(y_1, y_2)$ is polyhomogeneous as $r \to 0$, with $d_g(y_1, y_2) \sim r^{1-\alpha/2}$.
\end{corollary}
\begin{proof}
We have already shown that that there is a well-defined diffeomorphism which maps $(y_1,y_2)$ to $(x_0, y_0, \eta_0)$. We also analyzed that this
map has a smooth limit as $y_2 \to y_1$, i.e., for $r \to 0$.  In particular, $x_0$ depends smoothly on $r$. Next, by Lemma 9, the elapsed times 
$T(x_0, y_0, \pm \eta_0)$ to descend either of the two halves of this geodesic toward $y_1$ and $y_2$ are polyhomogeneous as $x_0 \to 0$.  
Finally, the lengths of these two half-geodesics $\gamma_\pm$ are computed using the usual formula
\[
\ell(\gamma_\pm) = \int_0^{T(x_0, y_0, \pm\eta_0)} |\gamma_\pm'(t)| \, dt,
\]
which is smooth in $T$, and hence polyhomogeneous in $x_0$ and thus also in $r$ .  Since $T \sim c x_0^{1-\alpha/2}$, this is the behavior of $d_g(y_1, y_2)$ as well.
\end{proof}

We note that it is possible to arrive at essentially the same conclusion, at least at the level of an estimate of order of growth but without the expansion,
by a more elementary method.   

\begin{proposition}
\label{prop:lengths-boundary} 
Then there are uniform constants $0 < C_1 < C_2$ so that
\begin{equation}
C_1 d_g(y_1,y_2) \leq d_{h_0}(y_1,y_2)^{1-\frac{\alpha}{2}} \leq C_2 d_g(y_1, y_2).
\end{equation}
for all $y_1,y_2\in\partial M$. 
\end{proposition}

\begin{remark}
Observe that this result shows in yet a different way that boundary measurements determine $\alpha$.
\end{remark}

\begin{proof}[Proof of proposition~\ref{prop:lengths-boundary}]
\label{pfsec:prop:lengths-boundary}
 It is convenient here to use coordinates $(s,y)$ where
\begin{equation}
\label{eqn:g-normal-form}
g = ds^2 + s^{-\beta } (1-\alpha/2)^{-\beta} h(s,y,dy) \quad \text{and} \quad h(0,y,dy) = h_0.
\end{equation}
In fact, $s = x^{1-\alpha/2}/(1-\alpha/2)$ and $\beta = 2\alpha/(2-\alpha)$.

Suppose, as before, that $y_1,y_2 \in \partial M$ are sufficiently close to one another.
We work locally near the boundary in the coordinates $(s,y)$ so that the metric is of the form~\eqref{eqn:g-normal-form},
and in particular the distance of a point $(s,y)$ to the boundary is $s$. 
Let $\gamma$ be the unique interior geodesic connecting these boundary points, and set $k = \max_t d_g(\partial M,\gamma(t))$. 
In the following, $C_\alpha$ denotes various positive constants depending on $\alpha$ but not $y_1, y_2$. 

We now approximate $\gamma$ by a ``quasi-geodesic''.
Define curves $\gamma_1, \gamma_2$ by $\gamma_j(t) = (t,y_j)$, $j = 1, 2$, 
$0 \leq t \leq \eps$, where $\eps$ is to be determined, and let $\gamma_3$ be the interior geodesic which connects $\gamma_1(\eps)$ to $\gamma_2(\eps)$. Denote by
$\gamma_\eps$ the concatenation of these three curves; this connects $y_1$ to $y_2$. Hence $d_g(y_1,y_2) \le \min_{\eps} \ell_g(\gamma_\eps)$,
the length with respect to $g$ of this piecewise curve. Furthermore, 
\begin{equation}
\label{eqn:lengths-bnd1}
l_g(\gamma_\eps) = l_g(\gamma_1) + l_g(\gamma_2) + l_g(\gamma_3) = 2\eps + C_\alpha\eps^{-\frac{\alpha}{2-\alpha}}d_h(y_1,y_2).
\end{equation}
The right-hand side of~\eqref{eqn:lengths-bnd1} is minimized at $\eps = C_\alpha d_h(y_1,y_2)^{1-\frac{\alpha}{2}}$. This gives
\begin{equation}
\begin{split}
d_g(y_1,y_2) &\le C_\alpha d_{h_0}(y_1,y_2)^{1-\frac{\alpha}{2}} + C_\alpha d_{h_0}(y_1,y_2)(d_h(y_1,y_2)^{1-\frac{\alpha}{2}})^{-\frac{\alpha}{2-\alpha}} \\
&= C_\alpha d_{h_0}(y_1,y_2)^{1-\frac{\alpha}{2}}.
\end{split}
\end{equation}

Next, choose $t_0$ so that $k = d_g(\partial M,\gamma(t_0))$. Writing $\gamma(t) = (s(t),y(t))$, then 
\begin{equation}
\label{eqn:lengths-bnd2}
\ell_g(\gamma) \ge \int_0^{\ell_g(\gamma)} \nabs{\dot s(t)} \,dt = \int_0^{t_0} \dot s(t) \,dt - \int_{t_0}^{\ell_g(\gamma)} \dot s(t) \,dt = 2k.
\end{equation}
On the other hand, since $s(t) \leq s(t_0)$ on the entire geodesic, 
\begin{equation}
\label{eqn:lengths-bnd3}
\begin{split}
\ell_g(\gamma)&\ge \int_0^{\ell_g(\gamma)} s(t)^{-\frac{\alpha}{2-\alpha}} \nabs{\dot y(t)}_{h_0} \,dt \ge s(t_0)^{-\frac{\alpha}{2-\alpha}}
\int_0^{\ell_g(\gamma)} \nabs{\dot y(t)}_{h_0} \,dt \\ &\ge k^{-\frac{\alpha}{2-\alpha}} d_{h_0}(y_1,y_2).
\end{split}
\end{equation}

Combining~\eqref{eqn:lengths-bnd2} and~\eqref{eqn:lengths-bnd3}, we obtain that $d_g(y_1,y_2) \ge C_\alpha d_{h_0}(y_1,y_2)^{1-\frac{\alpha}{2}}$,
while combining~\eqref{eqn:lengths-bnd1} and~\eqref{eqn:lengths-bnd3} yields $k \ge C_\alpha d_g(y_1,y_2)$. Therefore all $d_g$, $d_{h_0}$ 
and $k$ are comparable with constants only depending on $\alpha$ as claimed. 
\end{proof}

\begin{proposition}
\label{prop:bnd-hausdorff-dimension}
The Hausdorff dimension of $(\partial M, d_g)$ equals $\frac{2}{2-\alpha}(n-1)$. The Hausdorff dimension of $M$ equipped with this
same metric equals $\max\{n,\frac{2}{2-\alpha}(n-1)\}$.
\end{proposition}
\begin{proof}
By Proposition~\ref{prop:lengths-boundary}, $(\partial M, d_g)$ and $(\del M, d_{h_0}^{1-\frac{\alpha}{2}})$ are bi-Lipschitz equivalent, and
hence have the same Hausdorff dimension. It is enough then to compute $\mathrm{dim}\,(\partial M, d_{h_0}^{1-\frac{\alpha}{2}})$.
For simplicity, write this metric space as $(\del M, d_\alpha)$. 

It follows from the definition of Hausdorff measure that for all $\delta>0$, $\mathcal H^\delta(\del M, d_\alpha) = 
\mathcal H^{(1-\alpha/2)\delta}(\del M, d_{h_0})$, hence 
\begin{equation}
\dim_{\mathcal H} (\del M, d_\alpha) = \frac{2}{2-\alpha} \dim_{\mathcal H} (\del M, d_{h_0}) = \frac{2}{2-\alpha}(n-1).
\end{equation}
This proves the first claim. As for the second, this follows since $M = \partial M \cup \iM$ and $\dim_{\mathcal H} (\iM, d_g) = n$.
\end{proof}

\subsection{A travel time inverse problem for a gas giant}

As a first application of our study of geodesics on gas giants, we consider a preliminary inverse problem which asks whether the interior 
geometry of a gas giant can be recovered (modulo isometries) from knowledge of the Riemannian distances from interior points to the boundary. 
This simple application can be seen as a proof of concept of our gas giant geometry, leading to a proof as simple as that in the case $\alpha=0$.

The corresponding result is known both for compact Riemannian manifolds with boundary~\cite{KKL2001}
and in the Finsler setting~\cite{dHILS2019}. 
There is a more straightforward proof~\cite{ILS2023} in the Riemannian case when the metric is simple using a version of the 
Myers--Steenrod theorem from~\cite{dHILS2023}.
The result here is related to this simpler version. 

\begin{theorem}
\label{thm:bdy-distane-function}
Let $M$ be a compact manifold with boundary and, for $i = 1, 2$, suppose that $g_i$ are simple $\alpha_i$-gas giant metrics on $M$. Denote by
$d_i\colon M\times M \to \mathbb R^+$ the associated Riemannian distance functions. Define the maps $r_i\colon M \to \calC(\del M)$, where 
$r_i(x)$ is the function which sends $\del M \ni z \mapsto d_i(x,z)$. 

If the ranges of the two maps $r_1$ and $r_2$ are the same in $\calC(\del M)$, then $\alpha_1 = \alpha_2$ and $g_1$ is isometric to $g_2$ by
a diffeomorphism which is the identity on~$\del M$. 
\end{theorem}
\begin{proof}
First note that each map $r_i$ is well-defined, i.e., $r_i(x)$ is indeed a continuous function on $\del M$.
For standard incomplete metrics, 
this follows immediately from the triangle inequality.
In this setting, the same conclusion holds because, if $z, z' \in \del M$ and
$d_{h_i}$ is the distance function associated to the metric $h_i$ on $\del M$, then using the analysis of the last section gives the continuity estimate 
$|d_i(x,z) - d_i(x, z')| \leq d_{h_i}(z, z')^{1-\alpha_i/2}$.
We also observe the unique continuous extension of $d_i$ to the closed manifold
with boundary is also well-defined and continuous.

Next, it is also straightforward to check that each $r_i$ is injective. Indeed, if there were to exist two distinct points $x, x' \in M^\circ$ 
such that $r_i(x) = r_i(x')$, i.e., $d_i(x, z) = d_i(x', z)$ for all $z \in \del M$, then consider the maximally extended geodesic $\gamma$ 
which is length minimizing between any two of its points (this is where we use simplicity of the metrics) passing through $x$ and $x'$.
Suppose that one end of $\gamma$ meets $\del M$ at a point $z$, with $x'$ between $x$ and $z$. Then clearly $r_i(x)(z) = r_i(x')(z) + d_i(x,x')$,
so $r_i(x) \neq r_i(x')$ in $\calC(\del M)$.  The extended maps from all of $M$ are also injective. 

Continuing this same line of reasoning, we claim that in fact, if $x, x' \in M^\circ$, then $\norm{r_i(x)-r_i(x')}_\infty=d_i(x,x')$. 
This follows since, on the one hand, by the triangle inequality, $\norm{r(x)-r(x')}_\infty\leq d(x,x')$, while on the other, choosing
the minimizing geodesic $\gamma$ as above, then $d(z,x) - d(z,x') = \pm d(x,x')$, whence $\norm{r(x)-r(x')}_\infty\geq d_g(x,x')$.

As continuous injective maps from the compact Hausdorff space $(M, d_i)$ to $\calC(\del M)$, each $r_i$ is a homeomorphism onto
the common image $r_1(M) = r_2(M)$.  We may then define $\Psi\coloneqq r^{-1}_2 \circ r_1 \colon M \to M$.  By construction,
this is a bijective metric isometry. 

If $x \in \del M$, then $0 = r_1(x)(x) = r_2(\Psi(x))(x)$, so $x = \Psi(x)$, i.e., $\Psi$ is the identity on the boundary.  The fact that
$\Psi$ is a Riemannian isometry from $(M^\circ, g_1)$ to $(M^\circ, g_2)$ is then 
a consequence of the Myers--Steenrod theorem~\cite{MS1939,Pal2957}. 
\end{proof}

\section{Geodesic X-ray tomography on a gas giant}
This section studies the problem of unique reconstructibility of a function on a gas giant from the knowledge of its X-ray data i.e. integrals over 
all maximal geodesics. We prove that the X-ray data uniquely determines functions smooth up to the boundary.

The study of geodesic X-ray tomography in standard smooth Riemannian geometry originated in the work of Mukhometov~\cite{Muk1977}, who first 
proved the case $\alpha=0$ of our theorem~\ref{thm:injectivity-in-simple-geometry} below. For a comprehensive survey of the results, history 
and motivation of geodesic X-ray tomography, see~\cite{Sha1994,IM2019,PSU2023}. The X-ray transform is known to be injective on Cartan--Hadamard 
manifolds (see~\cite{LRS2018}) and in asymptotically hyperbolic geometry (see~\cite{GGSU2019}). These results are the closest relatives to our Theorem~\ref{thm:injectivity-in-simple-geometry}.

\begin{theorem}
\label{thm:injectivity-in-simple-geometry}
Let $M$ be a smooth manifold with boundary of dimension $n+1\geq2$.
Let $g = x^{-\alpha} \overline{g}$ be a gas giant metric on $M$ for some $\alpha \in (0,2)$ which is simple, i.e., non-trapping and free of
conjugate points. Suppose that a function $f \in \calC^\infty(\bar{M})$ 
has zero integral over all maximally extended $g$-geodesics. Then $f = 0$.
\end{theorem}

The proof of this theorem 
is based on a Pestov identity method. We begin by recalling relevant terminology, and refer to~\cite{Pat1999} for more details about the 
geometry of unit sphere bundles. 

\subsection{Pestov identity with boundary terms on a regular boundary}
Let $(M,g)$ be any compact smooth Riemannian manifold with smooth boundary (in this subsection $g$ is assumed to be smooth up to 
$\partial M$), and $S^\ast M$ its unit cosphere bundle.  This has the standard projection $\pi\colon S^\ast M \to M$, as well as a connection map 
$K \colon TS^\ast M \to TM$, defined by $K(\theta) = D_tc^\sharp(0)$; here $c$ is any curve in $S^\ast M$ with $c(0) = (x,\xi)$ and $\dot c(0) 
= \theta$, $c^\sharp(t) = c(t)^\sharp$ is the family of dual covectors, and $D_t$ is the Levi-Civita connection along $\pi(c(t))$.

There is an orthogonal decomposition 
\begin{equation}
TS^\ast M = \R X \oplus \mathcal H \oplus \mathcal V, \qquad 
\mathcal V = \Ker \der\pi, \quad \mathcal H = \Ker K.
\label{splitting}
\end{equation}
We denote by $N \to S^\ast M$ the bundle whose fibers are $N_{x,\xi} = \Ker \xi_x \subseteq T_xM$. The maps $\der\pi|_{\mathcal H} \colon \mathcal H \to N$ 
and $K|_{\mathcal V} \colon \mathcal V \to N$ are isomorphisms and we freely identify $\mathcal H \oplus \mathcal V = N \oplus N$.
We define the Sasaki metric $G$ on $S^\ast M$ by
\begin{equation}
G(\theta,\theta')
=
g(\der\pi(\theta),\der\pi(\theta'))
+
g(K(\theta),K(\theta'))
\end{equation}
for $\theta,\theta' \in T_{x,\xi}S^\ast M$. The splitting~\eqref{splitting} of $TS^\ast M$ is orthogonal with respect to~$G$.

The $G$-gradient of a smooth function $u$ on $S^\ast M$ can be written as
\begin{equation}
\nabla_Gu = (Xu,\gradh u,\gradv u)
\end{equation}
where the horizontal and vertical gradients $\gradh u$ and $\gradh u$ are smooth sections of the bundle $N$ and $X$ is the Hamiltonian vector field on $S^\ast M$.
The Riemannian curvature tensor maps sections $W$ of $N$ to sections of $N$ by the action
\begin{equation}
RW(x,\xi)
=
R(W(x,\xi),\xi^\sharp)\xi^\sharp.
\end{equation}
Let $\der\Sigma$ be the volume form of the Sasaki metric.

Now pull the volume form $\der\Sigma$ by the inclusion $\partial S^\ast M \to S^\ast M$ to get a natural volume form 
$\der\sigma$ on $\partial S^\ast M$. For all $u \in C^\infty(S^\ast M)$ define 
\begin{equation}
B(u) \coloneqq \int_{\partial S^\ast M} \nip{\gradv u}{\gradh u} + nuXu \,d\sigma.
\end{equation}
The following Pestov identity was proved in~\cite[p. 60]{GGSU2019}.

\begin{lemma}
\label{lma:pestov-with-bnd-term}
With notation as above, then for all $u \in C^\infty(S^\ast M)$, 
\begin{equation}
\label{eqn:pestov}
\nnorm{\gradv Xu}^2 = \nnorm{X\gradv u}^2 - \niip{R\gradv u}{\gradv u} + n \nnorm{Xu}^2 +  B(u). 
\end{equation}
\end{lemma}

\begin{remark}
By approximation, the identity~\eqref{eqn:pestov} continues to hold if $u \in C^1(S^\ast M)$ has $\gradv Xu \in L^2(N)$ and $X\gradv u \in L^2(N)$.
\end{remark}

\subsection{Proof of theorem~\ref{thm:injectivity-in-simple-geometry}}
\label{sec:proof-of-injectivity}

We return to the case where $g$ is a simple gas giant metric and present the proof of theorem~\ref{thm:injectivity-in-simple-geometry} using a collection 
of lemmas, the proofs of which appear in sections~\ref{sec:boundary-determination-xrt},~\ref{sec:derivates-of-u} and~\ref{sec:proof-of-pestov}.


\begin{lemma}
\label{lma:xrt-bnd-determination}
Let $g$ be a simple $\alpha$-gas giant metric on $M$, 
and let $f \in \calC^\infty(\bar{M})$. If $f$ integrates to zero over all maximally extended $g$-geodesics of $M$ then $f \in x^\infty \calC^\infty(M)$.
\end{lemma}

\begin{lemma}
\label{lma:regularity-of-u}
Let $g$ be a simple $\alpha$-gas giant metric on $M$, and suppose that $f \in x^\infty \calC^\infty(M)$. Then there is a solution $u \in x^\infty 
\calC^\infty(S^\ast\iM)$ to the transport equation $Xu=-f$ in $S^\ast\iM$ with $\nabla_G u \in x^\infty L^\infty(S^\ast M;TS^\ast M)$ and 
$X\gradv u, \gradv Xu \in L^2(N)$.
\end{lemma}

\begin{lemma}
\label{lma:pestov}
Let $u \in x^\infty C^\infty(S^\ast\iM)$ satisfy $Xu = -f$ as well as $\gradv Xu,X\gradv u \in L^2(N)$ and $\nabla_G u \in x^\infty L^\infty(S^\ast M;TS^\ast M)$. 
Then 
\begin{equation}
\nnorm{\gradv Xu}^2 = \nnorm{X\gradv u}^2 - \niip{R\gradv u}{\gradv u} + n \nnorm{Xu}^2.
\end{equation}
\end{lemma}

In the following, $\calC^\infty(N^\circ)$ denotes the space of sections of $N$ which are smooth over the interior $S^\ast\iM$. 

\begin{lemma}
\label{lma:global-index-form}
Let $g$ be a simple $\alpha$-gas giant metric on $M$. Then 
\begin{equation}
Q(W) = \nnorm{XW}^2 - \niip{RW}{W} \ge 0
\end{equation}
for all $W \in x^\infty \calC^\infty(N^\circ)$ with $W \in x^\infty L^\infty(S^\ast M;TS^\ast M)$.
\end{lemma}

\begin{proof}[Proof of theorem~\ref{thm:injectivity-in-simple-geometry}]
Suppose that $f \in \calC^\infty(M)$ integrates to zero over all maximally extended geodesics. Then $f \in x^\infty \calC^\infty(M)$ by 
Lemma~\ref{lma:xrt-bnd-determination} and so by Lemma~\ref{lma:regularity-of-u} there is a solution $u \in x^\infty \calC^\infty(S^\ast\iM)$ 
to the transport equation $Xu = -f$ with $\nabla_Gu \in x^\infty L^\infty(S^\ast M;TS^\ast M)$ and $\gradv Xu,X\gradv u \in L^2(N)$.
Apply the Pestov identity in Lemma~\ref{lma:pestov} to $u$ to get 
\begin{equation}
\label{eqn:pestov-in-action}
\nnorm{\gradv f}^2
=
Q(\gradv u)
+
n\nnorm{f}^2.
\end{equation}
Since $f$ is the lift of a function on $M$ to $S^\ast M$, $\gradv f \equiv 0$. In addition, $Q(\gradv u) \ge 0$. Thus by Lemma~\ref{lma:global-index-form}, 
the Pestov identity~\eqref{eqn:pestov-in-action} reduces to $0 \ge n\nnorm{f}^2$, so $f \equiv 0$ as claimed.
\end{proof}

\subsection{Boundary determination}
\label{sec:boundary-determination-xrt}

In this section we prove that a function smooth up to the boundary of $M$ is uniquely determined to any order at the boundary by its integrals over 
all maximal $g$-geodesics in $M$.   We first prove an auxiliary result about geodesics converging to a given boundary point.

\begin{lemma}
\label{lma:geodesics-to-a-point}
Let $g$ be a simple $\alpha$-gas giant metric on $M$. For any $\bar y \in \partial M$, there exists a sequence $\zeta_k \in S^*_{\bar y} M$
such that the lengths $l_g(\gamma_k)$ of the bicharacteristics $\gamma_k(t) = (z_k(t),\zeta_k(t))$ with $\gamma_k(0) = (\bar y, \zeta_k)$ are positive 
for all $k$ and converge to zero as $k \to \infty$.
\end{lemma}

\begin{proof}
Choose a smooth boundary curve $c \colon (-\eps,\eps) \to \partial M$ with $c(0) = \bar y$, and set $\bar y_k \coloneqq c(1/k)$. By simplicity, there 
is a unique unit speed bicharacteristic $\gamma_k(t) = (z_k(t),\zeta_k(t))$ with $z_k(0) = \bar y$, $z_k(\tau_k) = \bar y_k$; here $\tau_k$ is the exit time 
of $z_k$ (this is finite by lemma~\ref{lma:finite-exit-time}). We then let $\zeta_k \zeta_k(0)$. Since $\bar y_k \ne \bar y$, each $l_g(\gamma_k)$ 
has positive length. Moreover, by Proposition~\ref{prop:lengths-boundary} we have 
\begin{equation}
l_g(\gamma_k) = d_g^M(\bar y_k,\bar y) \le Cd_h^{\partial M}(\bar y_k,\bar y)^{1-\alpha/2},
\end{equation}
so the lengths converge to zero, as needed.
\end{proof}

We can now prove the following boundary determination lemma. We use arguments similar to the proof of~\cite[Theorem 2.1]{LSU2003}.
The only step in its proof where simplicity of the metric is needed is when Lemma~\ref{lma:geodesics-to-a-point} is invoked. 


\begin{proof}[Proof of lemma~\ref{lma:xrt-bnd-determination}]
We prove by induction on $\ell$ that for all $\bar y\in \partial M$ and every $\ell \geq 0$, $(\del_x^\ell f)(0,\bar y) = 0$. 

When $\ell = 0$, choose a sequence $\zeta_k \in S^*_{\bar y}M$ so that that the corresponding bicharacteristics $\gamma_k(t) = (z_k(t),\zeta_k(t))$
have lengths $l_g(\gamma_k)$ tending to zero, as in Lemma~\ref{lma:geodesics-to-a-point}. By hypothesis,
\begin{equation}
\frac{1}{\tau_k}\int_0^{\tau_k} f(z_k(t)) \,dt = 0,
\end{equation}
where $\tau_k$ is the length of the geodesic $z_k$.  Since $f$ is smooth, there exist $t_k \in (0,\tau_k)$ such that $f(z_k(t_k)) = 0$.  Clearly
$t_k < \tau_k \to 0$. Thus 
\begin{equation}
f(0,\bar y) = \lim_{k \to \infty} f(z_k(t_k)) = 0,
\end{equation}
as claimed.

Now assume, for any $\ell > 0$, that $\partial_x^jf(x,y)|_{(0,\bar y)} = 0$ for all $0 \le j < \ell$. We prove that $(\partial^\ell_xf)(0,\bar y) = 0$ by 
assuming the contrary, that $(\del_x^\ell f)(0,\bar y) \ne 0$ and arriving at a contradiction.

Assume that $(\partial_x^\ell f)(0,\bar y) > 0$. Since $f$ is smooth, $(\del_x^\ell f)(x,y) > 0$ for all $(x,y)$ in some neighborhood $\calU$ of $(0,\bar y)$.
Taking the Taylor expansion of $f$ at any $(0,y)$ and using the inductive hypothesis, we have that
\begin{equation}
\label{eqn:taylor-in-induction}
f(x,y) = x^\ell\del_x^\ell f(0,y) + \Order(x^{\ell+1}).
\end{equation}
By the positivity of the $\ell^{\text{th}}$ derivatives, there is a smaller neighbourhood $\bar y \in \calU' \subseteq \calU$ such that $f(x,y) > 0$ 
in $\calU'$.  Since $l_g(\gamma_k) \to 0$, the entire geodesic $z_k$ lies in $\calU'$ when $k$ is large. Hence the integral of $f$ over $z_k$ cannot 
vanish, a contradiction.

This proves that $f$ vanishes to order $\ell$ along $\del M$, and since this is true for all $\ell > 0$, we are done.
\end{proof}

As a corollary of this Lemma, we prove that the transport equation $Xu = -f$ admits a solution which is smooth in $M^\circ$, and that this solution 
vanishes to all orders at $\del M$ if $f \in \calC^\infty(M)$ is in the kernel of the X-ray transform.

Given $f \in \calC^\infty(\bar M)$, we define $u^f$ to be the function on $S^\ast M$ defined by the formula\footnote{In this section, unlike above, we denote the exit time by $\tau$ to adhere with the common convection.}
\begin{equation}
\label{eqn:uf}
u^f(z,\zeta) = \int_0^{\tau(z,\zeta)} f(\phi_t(z,\zeta)) \,dt;
\end{equation}
here $f$ is identified with its pullback $\pi^\ast f$, and $\phi_t(z,\zeta)$ is the cogeodesic flow.

\begin{corollary}
\label{cor:transport-eqn}
Let $g$ be a simple $\alpha$-gas giant metric on $M$, and let $f \in \calC^\infty(\bar M)$. If the integral of $f$ over all maximal geodesics in $M$ is zero,
then $u^f$ solves the transport equation $Xu = -f$ in $S^\ast\iM$, and $u^f \in x^\infty \calC^\infty (S^\ast\iM)$.
\end{corollary}

\begin{proof}
Since $f$, $\phi_t$ and $\tau$ are all smooth (see Lemma~\ref{lma:exit-time}), clearly $u^f \in \calC^\infty(S^\ast\iM)$.

We prove that $u^f(z,\zeta) = \Order(x^\ell)$ for all $\ell > 0$, where the constant depend only on $\ell$. It suffices to prove this at 
$(z_0,\zeta_0) \in S^\ast\iM$, so that $x_0 \in (0, \eps)$ and $\xi_0 < 0$.
For positive $\xi_0$ the claim follows from this one and vanishing integrals of $f$ over maximal geodesics.

Let $\gamma(t) = (z(t),\zeta(t))$ be a bicharacteristic $\gamma(0) = (z_0, \zeta_0)$.
We have already shown that $f(x,y) = \Order(x^\ell)$
for any $\ell$.
Since $x(t)$ is strictly decreasing by lemma~\ref{lma:monotonous-towards-b} and $g$ is non-trapping, we have
\begin{equation}
\begin{split}
\abs{u^f(z_0,\zeta_0)} &\le \int_0^{\tau(z_0,\zeta_0)} \abs{f(\phi_t(z_0,\zeta_0)} \,dt \le C_k \int_0^{\tau(z,\zeta)} x(t)^\ell \,dt \\
&\le \tilde C_kx(0)^\ell = \tilde C_kx_0^\ell
\end{split}
\end{equation}
for all $\ell > 0$, hence $u^f \in x^\infty \calC^\infty(S^\ast\iM)$.

To prove that $u^f$ solves $Xu = -f$ in $S^\ast\iM$, we compute just as for the classical case of metrics smooth up to the boundary. The point
is simply that $X$ differentiates along the cogeodesic flow and $u^f$ is defined by integration along the orbits of this flow.
\end{proof}

\subsection{Derivatives of the integral function}
\label{sec:derivates-of-u}
We now prove lemma~\ref{lma:regularity-of-u}. This involves an estimate of the derivatives of $u^f$, where $f$ has vanishing X-ray transform.
The first step is to show that normal Jacobi fields cannot blow up at the boundary with respect to the metric $\overline{g} = x^\alpha g$. 

\begin{lemma}
\label{lma:jacobi-growth}
Let $J(t)$ be a Jacobi field everywhere normal to a bicharacteristic curve $(z(t),\zeta(t))$ with $x(0) < \eps$ and $\xi(0) \le 0$.
Then $\abs{J(t)}_{\overline{g}} \le C$ and $\abs{D_tJ(t)}_{\overline{g}} \le Cx(t)^{-1}$ for all $t \in [0,\tau(z(0),\zeta(0))]$.
\end{lemma}
\begin{proof}
Choose any local coordinate system on $M$ near the endpoint of the projected geodesic. The Jacobi equation takes the form
\begin{equation}
\label{eqn:jacobi}
\ddot J^i + 2\Gamma^i_{jk}\dot \gamma^j\dot J^k + (\partial_k\Gamma^i_{jl})\dot \gamma^j\dot \gamma^l J^k = 0,
\end{equation}
where $\ddot J$ and $\dot J$ denote the usual derivatives of the coordinates of $J$ with respect to $t$. 
The Christoffel symbols of the actual metric $g$ are 
\begin{equation}
\begin{split}
\Gamma^0_{00}  &= -\frac{\alpha}{2}x^{-1},  \quad \Gamma^0_{i0} = 0, \quad \Gamma^m_{00} = 0 \\
\Gamma^0_{ij} &= -\frac{\alpha}{2}x^{-1}h_{ij} + \frac{1}{2}h_{ij} \\
\Gamma^m_{i0} &= -\frac{\alpha}{2}x^{-1}\delta^m_i + \frac{1}{2}h^{mk}\partial_xh_{ki} \\
\Gamma^i_{jk} &= \frac{1}{2}h^{mk} (\partial_jh_{ki} + \partial_ih_{kj} - \partial_kh_{ij}) \coloneqq H^i_{jk},
\end{split}
\end{equation}
where $H^i_{jk}$ is defined by this last equality. When $J(t)$ is normal to $z(t)$, the Jacobi equation reduces to 
\begin{equation}
\label{eqn:jacobi2}
\begin{split}
0 &= \ddot J^i + 2\Gamma^i_{0k}\dot x\dot J^k + 2\Gamma^i_{jk}\dot y^j\dot J^k + 2(\partial_k\Gamma^i_{0l})\dot x \dot y^l J^k +
(\partial_k\Gamma^i_{jl})\dot y^j \dot y^l J^k \\ 
&= \ddot J^i - \alpha x^{-1}\dot x \dot J^i + (h^{il}\partial_xh_{lk} \dot x + 2H^i_{jk}\dot y^j)\dot J^k \\
&\quad+ (\partial_k(h^{ip}\partial_xh_{pl})\dot x\dot y^l + (\partial H^i_{jl})\dot y^j\dot y^l)J^k.
\end{split}
\end{equation}
The coefficient of the third term on the right, involving $\dot J^k$, is bounded for $x \geq 0$, and since $\dot y = \Order(x^\alpha)$,  
equation~\eqref{eqn:jacobi2} becomes 
\begin{equation}
\label{eqn:jacobi3}
\ddot J^i - \alpha x^{-1}\dot x \dot J^i + F^i_k\dot J^k + x^\alpha G^i_k J^k = 0
\end{equation}
for some bounded functions $F^i_k$ and $G^i_k$. 

Equation~\eqref{eqn:jacobi3} can be reduced to a non-singular equation by rescaling $\dot J$. Define $W_1(t) = J(t)$ and $W_2(t) = x(t)^{-\alpha}\dot J(t)$,
so that $\dot W_1 = x^\alpha W_2$. Substituting into~\eqref{eqn:jacobi3} gives
\begin{equation}
0 = \alpha x^{\alpha - 1}\dot x W_2^i + x^\alpha \dot W^i_2 - \alpha x^{\alpha - 1}\dot x W^i_2 + x^\alpha F^i_kW^k_2 + x^\alpha G^i_kW^k_1
\end{equation}
which reduces to $\dot W^i_2=-F^i_kW^K_2-G^i_kW^k_1$.  This shows that $W = (W_1,W_2)$ satisfies $\dot W = AW$ where 
\begin{equation}
A = \begin{pmatrix} 0 & x^\alpha I \\ -G & -F \end{pmatrix} 
\end{equation}
and $I$ is the identity matrix, and $F = (F^i_k)$ and $G = (G^i_k)$.

It suffice to prove boundedness of the Jacobi field in the Euclidean metric $e$ with respect to the $(x,y)$ coordinates. We compute
\begin{equation}
\partial_t\abs{W(t)}^2_e = 2\dot W(t) \cdot W(t) = 2A(t)W(t) \cdot W(t).
\end{equation}
Since $A$ is continuous up to $\del M$, and hence bounded, we get $\partial_t\abs{W(t)}^2_e\le C\abs{W(t)}^2_e$. By Grönwall's inequality,
$\abs{W(t)}^2_e \le C\abs{W(0)}^2_e$.   This proves that $\nabs{J(t)}^2_e \le C$ and $\nabs{\dot J(t)} \le Cx(t)^{2\alpha}$, and hence 
$\abs{D_tJ(t)}_e \le Cx(t)^{-1}$, as claimed.
\end{proof}

By Lemma~\ref{lma:jacobi-growth}, we can now estimate derivatives of $u^f$.

\begin{lemma}
\label{lma:derivatives-uf}
If $f \in x^\infty \calC^\infty (M)$, then $\nabla_G u^f(z,\zeta) = \Order(x^\ell)$ for any $\ell \geq 0$, where the constants are uniform 
in $(z,\zeta) \in S^\ast\iM$.
\end{lemma}

\begin{proof}
It suffices to prove that $\partial_\theta u^f(z,\zeta) = \Order(x^\ell)$ uniformly on $S^\ast\iM$, where $x < \eps $, $\xi \le 0$ and $\theta \in 
T_{z,\zeta}S^\ast\iM$ with $\theta \perp X$. For convenience, identify $f$ with its lift $\pi^\ast f$ to $S^\ast\iM$.  Choose a smooth curve $c(s)$ 
in $S^\ast\iM$ with $c(0) = (z,\zeta)$ and $\dot c(0) = \theta$. Then 
\begin{equation}
\label{eqn:der-uf}
\begin{split}
& \partial_\theta u^f(z,\zeta) = \frac{d}{ds} \int_0^{\tau(c(s))} f(\phi_t(c(s))) \,dt \bigg|_{s = 0} \\
&= f(\phi_{\tau(z,\zeta)}(z,\zeta))\frac{d}{ds}\tau(c(s))\bigg|_{s=0} + \int_0^{\tau(x,\zeta)} \frac{d}{ds} f(\phi_t(c(s))) \bigg|_{s=0} \,dt.
\end{split}
\end{equation}
Since $\tau$ is smooth in $S^\ast\iM$ and $f$ vanishes on $\del M$, the first term on the right here vanishes. The second interior term
is estimated using Jacobi fields. 

Let $J_\theta(t)$ be the Jacobi field along the geodesic $\pi(\phi_t(z,\zeta))$ with initial conditions $J_\theta(0) = d\pi(\theta)$ and $D_tJ_\theta(0) = K(\theta)$.
In the splitting of $TS^\ast M$, the differential $d\phi_t(\theta)$ splits into $J_\theta(t) = d\pi(d\phi_t(\theta))$ and $D_tJ_\theta(t) = K(d\phi_t(\theta))$.
Thus 
\begin{equation}
\begin{split}
\frac{d}{ds} f(\phi_t(c(s))) \bigg|_{s=0} &= d(\pi^\ast f) (\partial_s\phi_t(c(s))|_{s = 0}) \\ &= \pi^\ast(df)(J_\theta(t),D_tJ_\theta(t)) \\ &= df(J_\theta(t)).
\end{split}
\end{equation}
Now, both $f$ and $df$ are $\Order(x^\ell)$ for all $\ell$. Since $\theta \perp X$, $J_\theta$ is normal to this geodesic, Lemma~\ref{lma:jacobi-growth}
implies that $|J_\theta(t)|_{\overline{g}}$ remains bounded. Thus the integrand in the second term is bounded by $Cx(t)^\ell$. Since $x$ is strictly decreasing 
on $[0,\tau(z,\zeta)]$ (and the metric is non-trapping), this shows that $\partial_\theta u^f(z,\zeta) = \Order(x^\ell)$ for all $\ell$, as claimed.
\end{proof}

\begin{proof}[Proof of lemma~\ref{lma:regularity-of-u}]
Let $f \in x^\infty \calC^\infty(M)$ and set $u = u^f$. By Corollary~\ref{cor:transport-eqn} and Lemma~\ref{lma:derivatives-uf}, the integral function $u$ 
satisfies $Xu = -f$ in $S^\ast\iM$ and $u \in x^\infty \calC^\infty(S^\ast\iM)$ and $\nabla_Gu \in x^\infty L^\infty(S^\ast M;TS^\ast M)$.
It remains to prove that $\gradv Xu,X\gradv u \in L^2(N)$.

Since $u$ solves the transport equation and the lift of $f$ to $S^\ast M$ depends only on $x$, we see that $\gradv Xu = -\gradv f = 0$, which is in $L^2(N)$. 
Now use the commutator formula $[X,\gradv ] = -\gradh $, valid in $S^\ast\iM$, (cf.\ ~\cite[Appendix A]{PSU2015}) to see that
\begin{equation}
\label{eqn:x-grav-estimate}
\nnorm{X\gradv u}_{L^2} = \nnorm{\gradh u}_{L^2} \le \nnorm{\nabla_G u}_{L^2}.
\end{equation}
Since $\nabla_G u \in x^\infty L^2(S^\ast M;TS^\ast M) \subset \nabla_G u\in L^2(S^\ast M;TS^\ast M)$, \eqref{eqn:x-grav-estimate} gives that 
$X\gradv u \in L^2(N)$.
This proves all of the assertions.
\end{proof}

\subsection{Proof of the Pestov identity}
\label{sec:proof-of-pestov}

We complete this entire argument by proving Lemmas~\ref{lma:pestov} and~\ref{lma:global-index-form}.

\begin{proof}[Proof of lemma~\ref{lma:pestov}]
Let $u \in x^\infty \calC^\infty(S^\ast\iM)$ be such that $\gradv Xu,X\gradv u \in L^2(N)$ and $\nabla_G u \in x^\infty L^\infty(S^\ast M;TS^\ast M)$. 
In adapted coordinates $(x,y)$ near $\del M$, consider the truncated manifold $M_\eps \coloneqq \{x \ge \eps\}$. The restriction of $g$ to
this truncation is smooth and non-degenerate up to $\del M_\eps$. By Lemma~\ref{lma:pestov-with-bnd-term}, for any $w \in \calC^\infty(S^\ast M_\eps)$,
\begin{equation}
\nnorm{\gradv Xw}^2_\eps = \nnorm{X\gradv w}^2_\eps - \niip{R\gradv w}{\gradv w}_\eps + (n-1) \nnorm{Xw}^2_\eps + B_\eps(w).
\end{equation}
In particular, this holds for the restriction of $u$ to $S^\ast M_\eps$. We prove that the identity on all of $S^\ast M$ by taking the limit $\eps \to 0$.

First, since $\nabla_G u \in x^\infty L^\infty(S^\ast M;TS^\ast M)$, we see that
\begin{equation}
\nabs{B_\eps(u)} \le C\eps^\ell\Vol(\{x = \eps\})
\end{equation}
for large $\ell$ in the sense of the inherited volume form of the submanifold $\{x=\eps\}$.
The volume of $(M,g)$ is finite when $\alpha < 2/n$;  if $\alpha = 2/n$, the volume of $M_\eps$ is asymptotic to $-C\log(\eps)$,
while for $\alpha > 2/n$ it is asymptotic to $C\eps^{1-n\alpha/2}$. Choose $\ell$ large, it is clear that $B_\eps(u) \to 0$ as $\eps \to 0$.

We next prove that the term involving curvature converges to the corresponding term in $S^\ast M$. The sectional curvatures of $g$ are asymptotic to 
$C_\alpha s^{-2}$ where $s$ is the distance to the boundary with respect to $g$; it is related to $x$ by $s = (1-\alpha/2)x^{1-\alpha/2}$. The pointwise 
inner product $\nip{R\gradv u}{\gradv u}$ is thus bounded by a multiple of $x^{-2+\alpha}\nabs{\nabla_G u}^2$. It follows that 
\begin{equation}
\niip{R\gradv u}{\gradv u} - \niip{R\gradv u}{\gradv u}_\eps \le C \int_{S^\ast M \setminus S^\ast M_\eps} x^{\ell-2+\alpha} \,d\Sigma 
\end{equation}
for $\ell$ large. We then compute that 
\begin{equation}
\begin{split}
\int_{S^\ast M \setminus S^\ast M_\eps} x^{\ell-2+\alpha} \,d\Sigma &= C \int_{M \setminus M_\eps} x^{\ell-2+\alpha-\frac{n\alpha}{2}} \,dxdV_h \\
&\le C \eps^{\ell-2+\alpha-\frac{n\alpha}{2}} \Vol_{\overline{g}}(M \setminus M_\eps). 
\end{split}
\end{equation}
Since $\ell$ can be chosen as large as desired, this last term vanishes as $\eps \to 0$, proving that
\begin{equation}
\niip{R\gradv u}{\gradv u}_\eps \to \niip{R\gradv u}{\gradv u}.
\end{equation}

Finally, since the pointwise norms $\nabs{X\gradv u}$, $\nabs{\gradv u}$ and $\nabs{Xu}$ are bounded by $\nabs{\nabla_Gu}$, a similar computation
shows that we can take limits in the remaining terms $\nnorm{X\gradv u}^2_\eps$, $\nnorm{\gradv Xu}^2_\eps$ and $\nnorm{Xu}^2_\eps$.
\end{proof}

\begin{proof}[Proof of lemma~\ref{lma:global-index-form}]
We prove finally that 
\begin{equation}
Q(W) = \nnorm{XW}^2 -\niip{RW}{W} \ge 0
\end{equation}
for all $W \in x^\infty \calC^\infty(N^\circ) \cap x^\infty L^\infty(S^\ast M;TS^\ast M)$.

Choose $\chi \in \calC^\infty([0,\infty))$ with $\chi = 1$ in $[2,\infty)$, $\chi = 0$ in $[0,1]$ and $0 \le \chi \le 1$ everywhere. We use the special adapted 
coordinates $(x,y)$. For $(z,\zeta) \in S^\ast M$, write $\chi_\eps(z,\zeta) = \chi(x/\eps)$, and define $W_\eps = \chi_\eps W$. Then $W_\eps$ is 
smooth in the interior of $S^\ast M$ and supported in $S^\ast M_{\eps} = \{(z,\zeta) \in S^\ast M\,:\, x \ge \eps\}$. We claim that $Q(W_\eps) \to Q(W)$ 
as $\eps \to 0$.

By the product rule,
\begin{equation}
\label{eqn:q-eps-formula}
\begin{split}
Q(W_\eps) &= \int_{S^\ast M_\eps} \abs{\chi_\eps}^2 \left( \nabs{XW}^2 - \nip{RW}{W} \right) \,d\Sigma_\eps\\
&\quad+ \int_{S^\ast M_\eps} \abs{(X\chi_\eps)W}^2 - 2\chi_\eps(X\chi_\eps) \nip{RW}{W} \,d\Sigma_\eps.
\end{split}
\end{equation}
The first term on the right converges to $Q(W)$ as $\eps \to 0$ by dominated convergence. It suffices to prove that the second term also vanishes 
as $\eps \to 0$.

The derivative $X\chi_\eps$ is supported in $\{\eps \le x \le 2\eps\}$. In addition, for all $(z,\zeta) \in S^\ast\iM$
\begin{equation}
X\chi_\eps(z,\zeta) = \frac{d}{dt} \chi(x(t)/\eps) \Bigg|_{t = 0} = \eps^{-1}\dot x(0)\chi'(x/\eps),
\end{equation}
so $\abs{X\chi_\eps} \le C\eps^{-1}$ and the integrand in the last term of~\eqref{eqn:q-eps-formula} is bounded by a constant multiple 
of $C\eps^{-1}(\abs{W}^2 + \ip{RW}{W})$ in $\{\eps \le x \le 2\eps\}$ and vanishes elsewhere.  Also, the sectional curvatures are asymptotic to 
$C_\alpha x^{-2+\alpha}$.  Since $W \in x^\infty L^\infty(SM;TS^\ast M)$, we can bound $C\eps^{-1}(\abs{W}^2 + \ip{RW}{W})$ by a multiple of 
$\eps^{\ell - 3 + \alpha}$ in $\{\eps \le x \le 2\eps\}$, and hence the second integral in~\eqref{eqn:q-eps-formula} is bounded by
\begin{equation}
C\eps^{\ell-3+\alpha}\Vol_G(\{\eps \le x \le 2\eps\}) = \tilde C\eps^{\ell-3+\alpha}\Vol_g(\{\eps \le x \le 2\eps\}).
\end{equation}
The volume grows no faster than a fixed power of $\eps$, so choosing $\ell$ sufficiently large, we see that this term also vanishes in the limit. 
Thus $Q(W_\eps) \to Q(W)$ as $\eps \to 0$.

The final step is to note that since the $W_\eps$ are smooth and compactly supported in $S^\ast M_\eps$, and since the truncated manifold $M_\eps$ is simple
in the traditional sense, it follows from~\cite[Lemma 11.2]{PSU2015} that $Q(W_\eps) \ge 0$ for all $\eps > 0$. Thus its limit $Q(W)$ is also non-negative.
\end{proof}

\section{The Laplacian of $g$}
\label{sec:4}

We now turn to the final major theme of this paper, which is to determine a few of the fundamental analytic properties of the scalar Laplace--Beltrami 
operator $\Delta_g$ associated to a gas giant metric. This operator degenerates at $x=0$, hence is poorly behaved from the point of view of 
classical theory. However, as we explain here, it can be regarded as an operator with a ``uniform degeneracy'' as $x \to 0$, and as such, can be transformed
to lie in a class of operators for which there is already an extensive theory. We describe this transformation of $\Delta_g$ into an ``elliptic $0$-differential operator'', 
as studied in~\cite{Maz1991} (and elsewhere). Quoting results from that theory, we study some of basic mapping and regularity properties of $\Delta_g$. 

We begin by deriving an expression for this operator in terms of the Laplacian of the metric $\overline{g}$.  First observe that 
\[
g^{ij} = x^\alpha \overline{g}^{\, ij}, \qquad \det (g_{ij}) = x^{-\alpha n} \det (\overline{g}_{ij}).
\]
For simplicity, write $\det g = \det (g_{ij})$ and $\det \gbar = \det (\overline{g}_{ij})$. Using the usual special adapted
coordinates $z = (z_0, z') = (x,y)$, we compute
\[
\Delta_g =  x^{\alpha n/2} \frac{1}{\det \gbar} \, \del_{z_i} \left( x^{\alpha(1-n/2)} (\det \gbar) \,\, \gbar^{ij} \,\del_{z_j} \right) = 
x^\alpha \Delta_{\gbar} +  x^{\alpha-1} \alpha(1-\frac{n}{2})\, \del_x 
\]
As noted earlier, this operator is clearly degenerate at $x=0$. 

We now set this into the context of the class of uniformly degenerate, or $0$-, differential operators. Using coordinates $(x,y)$
near $\del M$, we recall that a differential operator $L$ is called a $0$-operator if it can be expressed locally as a linear combination 
of products of smooth vector fields, each of which vanish at $\del M$.  The space of all smooth vector fields vanishing at $\del M$ is 
denoted $\mathcal V_0(M)$, and called the space of $0$ vector fields. It is generated over $\calC^\infty(M)$ by the `basis' vector fields 
$x\del_x, x\del_{y_1}, \ldots, x \del_{y_{n-1}}$, i.e., 
\[
\mathcal V_0(M) = \mbox{span}_{\calC^\infty} \, \{ x\del_x, x\del_{y_1}, \ldots, x\del_{y_{n-1}}\}.
\]
Thus a $0$-operator can be written in small neighborhoods as a finite sum of smooth multiples of products of these vector fields. In particular, 
for example, a $0$-operator of order $2$ is one which takes the form
\[
L = \sum_{ j + |\beta| \leq 2}  a_{j\beta } (x,y) (x\del_x)^j (x\del_y)^\beta.
\]

For the present purposes, the key point is that $\Delta_g$ assumes this form after multiplication by the factor $x^{2-\alpha}$.  (In carrying out some 
of the arguments below, it is occasionally more transparent to maintain symmetry of the operator by pre- and post-multiplying by $x^{1-\alpha/2}$, but
we shall not get into this level of detail). To illustrate this, let $x$ be a special boundary 
defining function, so that $\gbar = dx^2 + h$, where $h(x)$ is a smooth family of metrics on $\del M$, pulled back to this collar neighborhood 
by the projection $(x,y) \mapsto y$.  Then
\begin{equation}
\begin{aligned}
\Delta_g \coloneqq  & x^{\alpha} \Delta_{\bar{g}} - \alpha (n/2-1) x^{\alpha-1} \del_x \\ & = 
x^\alpha(\del_x^2 + q(x,y) \del_x + \Delta_{h(x)}) - x^{\alpha-1}\alpha( n/2 - 1) \del_x.
\end{aligned}
\end{equation}
where $q(x,y)$ is related to derivatives of $\det h$, but its precise expression is irrelevant since it is a higher order error term.  This operator is symmetric with 
respect to the measure $x^{-\alpha n/2} dx dV_{h(x)}$.  From this we see directly that $x^{2-\alpha}\Delta_g$ is a $0$-differential operator.

Associated to a $0$-differential operator is its $0$-symbol, which is obtained by writing $L$ as a sum of products of the generating vector fields, 
and then replacing each $x\del_x$ by $\xi$ and $x\del_{y_i}$ by $\eta_i$, and then dropping all terms with homogeneity less than that of the degree of $L$. 
In particular, for $L = x^{2-\alpha} \Delta_g$,
\[
{}^0\sigma_2(L)(x,y; \xi, \eta) \coloneqq  \xi^2 + |\eta|_{h(x)}^2.
\]
This is not apparently an invariant definition, but $(\xi,\eta)$ turn out to be natural linear variables on the fiber of a certain replacement for 
the cotangent bundle $T^*M$, and ${}^0\sigma_2(L)$ is a well-defined homogeneous polynomial of degree $2$ on these linear fibers. 
In any case, a $0$-operator is called $0$-elliptic if this symbol is non-vanishing (or invertible, if a system) when $(\xi,\eta)\neq (0,0)$; 
this operator $L$ is obviously $0$-elliptic.   As such,  the calculus of $0$-pseudodifferential operators offers analogues
of all the familiar constructions in pseudodifferential theory.  In particular, there is an elliptic parametrix construction for $L$, and
the various properties of the parametrix $G$ for $L$ obtained through this construction lead to sharp mapping and regularity properties which
are used below.

Now consider the densely defined unbounded operator 
\begin{equation}
\Delta_g\colon  L^2(M, dV_g) \longrightarrow L^2(M, dV_g).
\label{L2L2}
\end{equation}
This is symmetric on the core domain $\calC^\infty_0( M^\circ)$ of smooth functions compactly supported in the interior of $M$, and one of
the starting points for the analysis of the Laplacian is to determine whether this symmetric operator has a unique self-adjoint extension,
or if boundary conditions need to be imposed to obtain a self-adjoint realization?  Once that is accomplished, one may proceed to 
study the spectrum of any such self-adjoint extension.

We first recall some facts relevant to determining whether $\Delta_g$ is essentially self-adjoint.  In the following, we translate
some definitions from the development of the $0$-calculus to the present setting (rather than working directly with
the $0$-operator $L = x^{2-\alpha}\Delta_g$ simply to avoid too many confusing changes of notation. 

A fundamental invariant of $\Delta_g$ in this geometric setting is its pair of indicial roots, $\gamma_\pm$.  These are the values $\gamma$ 
such that solutions of $\Delta_g u$ grow or decay like $x^\gamma$. More formally, these are the exponents which yield approximate solutions 
in the sense that 
\[
\Delta_g x^\gamma = \mathcal O(x^{\gamma -1 + \alpha})
\]
rather than the expected rate $\mathcal O(x^{\gamma - 2 + \alpha})$. In other words, $\gamma$ is an indicial root if there is some leading order
cancellation.  To calculate these, we compute
\[
\Delta_g x^\gamma = \left( \gamma (\gamma-1) - \alpha(n/2 - 1) \gamma\right)x^{\gamma-2+\alpha} + \mathcal O(x^{\gamma -1 + \alpha}),
\]
and hence $\gamma$ must satisfy $\gamma^2 - ( \alpha(n/2-1) + 1) \gamma = 0$, or finally
\[
\begin{split}
\gamma_\pm   = &  0, \alpha(n/2-1) + 1 \\ 
= & \frac12( \alpha(n/2-1) + 1) \pm \frac12(\alpha(n/2-1) + 1).
\end{split}
\]
This last expression is included to emphasize the symmetry of $\gamma_\pm$ around their average, which is useful below. 

Next, observe that a function $x^\gamma$ lies in $L^2(dV_g)$ near $x=0$ if and only if 
\[
\gamma  > \frac12 ( n \alpha/2 - 1). 
\]
We call this threshold the ``$L^2$ cutoff weight''.   It is most natural to let $\Delta_g$ act on the Sobolev spaces adapted
to the $0$-vector fields:
\[
H^k_0(M, dV_g) = \{u:  V_1 \ldots V_\ell u \in L^2(dV_g)\ \ \ell \leq k,\ \ \mbox{each}\ V_i \in \mathcal V_0(M) \},
\]
and their weighted version $x^\mu H^k_0 = \{u = x^\mu v: v \in H^k_0\}$. It clear from this definition that 
\[
\Delta_g\colon  x^\mu H^2_0 \longrightarrow x^{\mu - 2 + \alpha} L^2
\]
is bounded for every $\mu$. In particular, $\Delta_g u \in L^2$ if $u \in x^\mu H^2_0$ where $\mu \geq 2-\alpha$. Since
$\calC^\infty_0(M^\circ)$ is dense in $x^{2-\alpha}H^2_0$, it is clear that the minimal domain, i.e., the minimal closed extension
from the core domain, of \eqref{L2L2} is contained in $x^{2-\alpha}H^2_0$.  Using the parametrix for $\Delta_g$ alluded to above,
it can be proved that this is an equality:
\[
\mathcal D_{\mathrm{min}}(\Delta_g) = x^{2-\alpha}H^2_0(M, dV_g).
\]
On the other hand, we also define the maximal domain $\mathcal D_{\mathrm{max}} = \{u \in L^2: \Delta_g u \in L^2\}$.

\begin{proposition}
The operator $\Delta_g$ is essentially self-adjoint on $L^2$, i.e., 
\[
\mathcal D_{\mathrm{min}} = \mathcal D_{\mathrm{max}}
\]
if and only if $\alpha > 2/n$.
\end{proposition}
\begin{proof}
They key issue is whether either of the indicial roots $\gamma_\pm$ lie in the critical weight interval
\[
\mu_- \coloneqq
\frac12 (n\alpha/2 - 1)  \leq \mu \leq \frac12 (n\alpha/2 - 1) + 2-\alpha
\eqqcolon \mu_+.
\]
Notice that the midpoint of this critical interval is $\frac12 (n\alpha/2 - 1) + 1-\alpha = \frac12 ( \alpha (n/2 - 1) + 1)$, which
is precisely the same as the midpoint of the gap between the two indicial roots.  The width of this weight interval is $2-\alpha$,
whereas $\gamma_+ - \gamma_- = \alpha (n/2 - 1) + 1$.  We claim that 
\[
\gamma_- < \mu_- < \mu_+ < \gamma_+
\]
precisely when $\alpha > 2/n$, which is verified by noting that $\alpha (n/2-1) + 1 > 2-\alpha$ precisely then.

The relevance of whether the indicial roots are included in the critical weight interval is that, using the parametrix carefully, one
can deduce that if $\gamma_\pm$ do not lie in this critical weight interval, then $u \in L^2$ and $\Delta_g u \in L^2$ imply
that $u \in x^{2-\alpha} H^2_0 = \mathcal D_{\mathrm{min}}$.  However, when $\alpha \leq 2/n$, then we can only deduce that
\[
u(x,y) \sim \sum a_j(y) x^{\gamma_- + j} + \sum b_j(y) x^{\gamma_+ + j}.
\]
This asymptotic expansion has some complicating features, such as that if $a_0 \not\equiv 0$, then the coefficients $a_j, b_j$ may
only have finite regularity (and will have negative Sobolev regularity for large $j$.  Conversely, there exists a solution of $\Delta_g u = 0$
where $u$ has an expansion of this type with any prescribed smooth leading coefficient $a_0(y)$. In any case, the upshot is that the maximal
domain is far bigger than the minimal domain in this case. 
\end{proof}

When $\alpha < 2/n$, there are many possible self-adjoint extensions. The most prominent, and the one we shall use below, is the Dirichlet
extension. This corresponds to the choice of domain $\mathcal D_{\mathrm{Dir}}$ consisting of those $u \in L^2$ such that $\Delta_g u \in L^2$
and where the leading coefficient $a_0(y)$ in the expansion above vanishes.   Other self-adjoint extensions correspond to other types
of conditions on the pair of leading coefficient $(a_0(y), b_0(y))$.  We do not detail these below, except mentioning the most standard
other ones: the Neumann extension, where $b_0(y) \equiv 0$, and the family of Robin extensions, corresponding to conditions of the
form $A(y) a_0(y) + B(y) b_0(y) \equiv 0$, where $A, B$ are given smooth functions.

\begin{proposition}
\label{prop:Fredholm-and-discrete-spectrum}
Let $\mathcal D$ be a domain of self-adjointness for $\Delta_g$ as above. Then $(\Delta_g, \mathcal D)$ is a Fredholm operator on $L^2$ with discrete spectrum. 
\end{proposition}
\begin{proof}
The first step is to show that this operator is Fredholm.  This follows from the existence of its parametrix. This is a $0$-pseudodifferential operator
of order $-2$ which maps $L^2$ onto $\mathcal D$ (possibly modulo compact errors), and which satisfies $G \circ \Delta_g = \mathrm{Id} - R_1$,
$\Delta_g \circ G = \mathrm{Id} - R_2$, where $R_1$ and $R_2$ are compact operators on $L^2$ and on $\mathcal D$ (with its graph
topology) respectively.   As noted earlier, the construction of this parametrix is one of the standard consequences of $0$-ellipticity; details are 
given in \cite{Maz1991}.  When $\alpha < 2/n$, a slightly more intricate construction is needed which incorporates the choice of boundary 
conditions; this appears in~\cite{Maz-Ver}.

The key point here is that the operator $G$ is constructed as an element of the $0$-pseudodifferential calculus. This means that its Schwartz 
kernel is a very well-understood object which, as a distribution on $M \times M$, has explicit asymptotic expansions at the boundary faces 
of this product, and a slightly more intricate, but equally explicit expansion near the corner of $M^2$.   The precise details are omitted.
The upshot, however, is that it then follows by general properties of such pseudodifferential operators proved in~\cite{Maz1991} that
$G$ is bounded on $L^2$. Of course, as a pseudoinverse to $\Delta_g$,  its range must lie in $\mathcal D$.  Since it is a ($0$-)pseudodifferential
operator of order $-2$, it is clear that the elements in $G(L^2)$ have two derivatives in $L^2$, at least in the interior of $M$. However, slightly
more is true, and the precise statement is that for any $f \in L^2$ and any two vector fields $V_1, V_2 \in \mathcal V_0(M)$,  we must 
have that $V_1 V_2 (G f) \in x^\eps L^2$ for some fixed $\eps > 0$ which is independent of $f$.  This is summarized by saying
that $G\colon L^2 \rightarrow x^\eps H^2_0$, where the range is a weighted $0$-Sobolev space.  We may then invoke the $L^2$ version
of the Arzel\`{a}--Ascoli theorem, which may be used to prove that $x^\eps H^2_0 \hookrightarrow L^2$ is a compact embedding. 
This shows that the domain of $\Delta_g$ is compactly contained within $L^2$, and hence that $(\Delta_g, \mathcal D)$ has
discrete spectrum. 

We say a few more words about this parametrix construction, particularly  when $\alpha > 2/n$. Write $\Delta_g = x^{\alpha/2 - 1} L x^{\alpha/2 - 1}$;
as noted earlier, $L$ is an elliptic $0$-operator.  The singular factor has been distributed on opposite sides of $L$ to preserve symmetry. 
Let $\overline{G}$ be a parametrix for $L$ as constructed in~\cite{Maz1991}. Thus $\mathrm{Id} - L \overline{G}  = R'_1$, $\mathrm{Id} - \overline{G} L = R_2'$, 
where $R_1'$, $R_2'$ are operators with smooth kernels on the interior of $M \times M$, and which admit classical expansions at all boundary
faces of a certain resolution (or blow-up) of this product, with coefficients in these expansion smooth functions on the corresponding
boundary faces.  We then write $G = x^{1-\alpha/2} \overline{G} x^{1-\alpha/2}$, so that $\Delta_g \circ G = \mathrm{Id} - x^{\alpha/2-1} R_1 x^{\alpha/2 - 1} = 
I - R_1$, $G \circ \Delta_g = \mathrm{Id} - x^{\alpha/2-1} R_2' x^{\alpha/2-1} = \mathrm{Id} - R_2$.  These remainder terms are much better,
inasmuch as they have smooth Schwartz kernels which have polyhomogeneous expansions at the two boundary hypersurfaces of $M^2$,
without need for the resolution (or blow-up) process.  

If $\Delta_g u = f \in L^2$, then applying $G$, we get that $u = R_1 u + G f = R_1 u + x^{1-\alpha/2} \overline{G} x^{1-\alpha/2} f$.  The first term is 
polyhomogeneous on $M$, and decays at a fixed rate strictly greater than the $L^2$ cutoff. When $\alpha > 2/n$, the range of $G$ lies in 
$x^{2-\alpha}H^2_0$. This range is identified with the domain of self-adjointness $\mathcal D$ (again, when $\alpha > 2/n$), hence, 
as described above, $\mathcal D \subset L^2$ is indeed compact. 
\end{proof}

We now take up our final problem. Fix a domain $\mathcal D \subset L^2$ where $(\Delta_g, \mathcal D)$ is self-adjoint.  To be very concrete below,
we assume that this is the Dirichlet extension henceforth. As just proved, the Dirichlet Laplacian has discrete spectrum $0 \leq \lambda_0 < \lambda_1 \leq 
\lambda_2  \leq \ldots$. 

Next consider the truncated manifold $M_\eps = \{ p \in M: x(p) \geq \eps\}$, where $x(p)$ is just the value of the boundary defining function
$x$ at $p$ (we assume that $x$ has been extended to be a smooth function on the interior of $M$ which is strictly positive on $M_\eps$). 
Then $\Delta_g$ restricts to an operator acting on $H^2(M_\eps)$ functions which vanish at $\del M_\eps$. This operator has discrete spectrum 
as well, by classical elliptic theory, and we denote its eigenvalues by $0 < \lambda_0(\eps) < \lambda_1(\eps) \leq \ldots$.  By classical perturbation theory, 
each $\lambda_j(\eps)$ can be regarded as a continuous, and piecewise smooth, function of $\eps$.  The basic question is whether the spectrum of 
$(\Delta_g, \mathcal D_{\mathrm{Dir}})$ on $M_\eps$ converges to the spectrum of $(\Delta_g, \mathcal D)$ on $M$.    While there are such statements that
can be made about the entire spectrum at once, we consider here the variation of individual eigenvalues.

\begin{proposition} \label{Prop:truncation}
For each $j = 0, 1, 2, \ldots$, the function $\lambda_j(\eps)$ converges to $\lambda_j$ as $\eps \to 0$. In fact, there exists a constant $C_j > 0$ such that
\[
|\lambda_j(\eps) - \lambda_j| \leq C_j \eps^{\alpha(n/2 - 1) + 1}
\]
\end{proposition}

\begin{proof}
We will denote $\partial_\eps$ by dot and $\partial_x$ by prime.

Let us focus on a particular eigenvalue $\lambda_j(\eps)$. For simplicity, we first make the computations below assuming that this is a simple eigenvalue, staying
away from eigenvalue crossings.    Thus, dropping the index $j$, assume that $\Delta_g \phi = \lambda(\eps) \phi$ on $M_\eps$, with $\phi = 0$ on $\del M_\eps$. 
We shall construct a family of diffeomorphisms $F_\eps\colon M \to M_\eps$, with $F_0 = \mathrm{Id}$. Using these to pull back all the data on $M_\eps$, we consider
the family of metrics $g_\eps = F_\eps^* g$, the associated Laplace operators $\Delta_{g_\eps}$, and eigenfunctions $\phi_\eps$ which are smooth functions on $M$
vanishing at $\del M$.  Choosing these to have $L^2(M, dV_\eps)$ norms equal to $1$, the proof involves estimating the quantity
\[
\dot{\lambda} = \int_M (\dot{\Delta}_{g_\eps} \phi_\eps) \phi_\eps\, dV_{\eps}.
\]

There are two parts to this. In the first, we obtain the uniform estimate 
\[
|\phi_\eps| \leq C x (\eps + x)^{\alpha(n/2 - 1)},
\]
with a constant $C$ independent of $\eps$. Thus $\phi_\eps$ vanishes only like $x$ when $\eps > 0$, but like $x^{\alpha(n/2-1)+1} = x^{\gamma_+}$ when $\eps = 0$. 
In the second, we must compute $\dot{\Delta}$.

To get started, we define the diffeomorphisms $F_\eps$.  Using a specially adapted boundary defining function, define $\mathcal U = \{x < c\}$ for some small $c > 0$,
and identify $\mathcal U$ with $[0,c) \times \del M$. Define $F_\eps(x) = x + \eps \chi(x/\eps)$, where $\chi(s)$ is a smooth monotone non-negative function which equal $1$ for $s \leq 1$
and $0$ for $s \geq 4$. We also require that $|\chi'(s)| \leq 1/2$ for all $s$.  We then have that
\[
g_\eps = x_\eps^{-\alpha} ( dx_\eps^2 + h(x_\eps)).
\]
Define
\[
\frac{dx_\eps}{dx} \coloneqq J = 1 + \chi'(x/\eps), \ \ \mbox{and}\ \ \  \dot{x_\eps} = \frac{dx_\eps}{d\eps} = - (x/\eps) \chi'(x/\eps).
\]
Then, since $\del_{x_\eps} = J^{-1} \del_x$, we obtain
\[
\Delta_\eps = x_\eps^\alpha J^{-2} \del_x - x_\eps^{\alpha-1} ( \alpha(n/2-1) J^{-1} - x_\eps J'/J^2) \del_x + x_\eps \Delta_h.
\]
We have actually made the tacit assumption here that $\bar{g} = dx^2 + h$ with $h$ independent of $x$. The extra terms which appear when $h$ depends on $x$
are lower order in all the computations below, so we can safely omit them. 

Computing further, we arrive at the expression
\[
\dot{\Delta} = x_\eps^{\alpha-1} A \del_x^2 + x_\eps^{\alpha-2} B \del_x  + x_\eps^{\alpha-1} C \Delta_h,
\]
where $A$, $B$ and $C$ are expressions which are sums of terms, each a smooth bounded multiples involving the quantities $x_\eps \dot{J}$, $\dot{x}_\eps$, $x_\eps J'$ and 
$(x_\eps J')\dot{}$.
The key point here is that each of these terms is uniformly bounded in $\eps$, and supported in the region $\eps \leq x \leq 4\eps$. 

Now, granting the uniform estimate on the $\phi_\eps$ stated above, we see that $\dot{\Delta} \phi \sim \eps^{ \alpha-1 + \alpha(n/2-1) -1}$, and as before, supported in
$x \in [\eps, 4\eps]$. Thus
\[
\begin{split}
\int (\dot{\Delta} \phi) \phi\, dV_\eps & \sim \int_\eps^{4\eps}  \eps^{\alpha-1 + \alpha(n/2-1) - 1 + \alpha(n/2-1) + 1 - n\alpha/2}\, dx \\ & = \int_\eps^{4\eps} \eps^{\alpha(n/2-1) - 1}\, dx 
= 4 \eps^{\alpha(n/2-1)},
\end{split}
\]
as claimed. 

It remains to verify the assertion about the uniform bound on $\phi_\eps$. 
We indicate the more elementary of the two arguments.   First note that since the $L^2$ norm of $\phi_\eps$ equals $1$, and interior estimates
bound $|\nabla \phi_\eps$ on any subset $\{x \geq c > 0\}$, the functions $\phi_\eps$ are uniformly bounded on any compact subset of
the interior of $M$, and all vanish at the boundary.  We obtain a uniform upper bound of the form $|\phi| \leq C x (\eps + x)^\beta$
for any $\beta < \alpha(n/2-1)$. Since we may take $\beta$ arbitrarily close to this upper limit, this suffices to give the eigenvalue variation limit
above with arbitrary small loss in the exponent.

Now suppose that there is no uniform constant $C$ such that $|\phi_\eps(x,y)| \leq C x (\eps + x)^{\beta}$.
The bound is clearly true for $\eps \geq \eps_0 > 0$ and for $x \geq c > 0$, so  there must exist sequences $(x_j, y_j)$ and $\eps_j$ with 
$x_j \to 0$, $\eps_j \to 0$, such that after multiplying by a sequence of factors $1/C_j \to 0$, and writing $\phi_j$ instead of $\phi_{\eps_j}$, we have
\[
|\phi_j(x,y)| \leq x (\eps_j + x)^{\beta}, \ \ |\phi_j(x_j, y_j)| = x_j (\eps_j + x_j)^{\beta}. 
\]
Now rescale, setting $s = x/x_j$, $w = (y-y_j)/x_j$, to write this as
\[
|\phi_j(x,y)| \leq  x_j s (\eps_j + x_j s)^{\beta}.
\]

We now separate into two cases.  In the first, $\eps_j \gg x_j$, so we rewrite the right hand side of this inequality as $x_j \eps_j^\beta s (1 + (x_j/\eps_j) s)^\beta$.
Replacing $\phi_j$ by $\tilde{\phi}_j = \phi_j/x_j \eps_j^{\beta}$, we see that
\[
|\tilde{\phi}_j(s,w)| \leq s (1 + (x_j/\eps_j) s)^\beta,
\]
with equality at $(1,0)$. Taking a limit as $j \to \infty$, we conclude the existence of a limit $\phi_\infty$ which satisfies $|\phi_\infty| \leq s$ for
$s \geq 0$ and all $w \in \RR^{n-1}$.  Each of the $\phi_j$ is smooth up to $\del M_{\eps_j}$, and there is a uniform bound on the tangential derivatives
(this follows from the parametrix methods); this implies that $\phi_\infty$ is in fact constant in $w$, and so must satisfy the ODE
$s^2 \del_s^2 \phi_\infty - \alpha(n/2-1) s \del_s \phi_\infty = 0$, whence $\phi_\infty(s) = C s^{1 + \alpha(n/2-1)}$. This contradicts the bound above.

The other case is when $x_j \geq C \eps_j$ as $\eps_j \to 0$. Now rewrite the right hand side of the inequality as $x_j^{1 + \beta} s ( (\eps_j/x_j) + s)^\beta$.
Now normalize by dividing by $x_j^{1+\beta}$ to define $\tilde{\phi}_j$, and take a limit as before. This yields a function $\phi_\infty$ such tat
\[
|\phi_\infty(s,w)| \leq  s (c + s)^\beta,
\]
with equality at $(1,0)$, and where $c$ is the limit of (some subsequence) of the $\eps_j/x_j$. This constant is finite, and possibly $0$.  As before, $\phi_\infty$
is independent of $w$ and must equal a constant times $s^{\alpha(n/2-1) + 1}$ for all $s \geq 0$, which is inconsistent with this limiting bound as $s$ gets large. 

As noted earlier, there is a more sophisticated way to obtain a sharper bound, and in fact complete asymptotic expansions for $\phi_\eps$ 
as both $x \to 0$ and $\eps \to 0$.  This requires a generalization of the parametrix machinery described above. This generalization allows one
to treat not only degenerate operators such as $\Delta_g$, but also families of {\it degnerating} operators $\Delta_{g_\eps}$.  However, for simplicity
we do not describe or develop this point of view here.   What we have proved with this more elementary argument is the slightly weaker estimate that
each eigenvalue $\lambda(\eps)$ satisfies
\[
|\dot{\lambda}(\eps)| \leq C_\delta \eps^{\alpha(n/2-1) + 1-\delta}
\]
for any $\delta > 0$. 

We then return to the case of a degenerate eigenvalue $\lambda$.
The eigenspace at $\eps=0$ has the orthonormal basis $\{\phi^1,\dots,\phi^m\}$ so that each $\phi^k$ is the limit of eigenfunctions $\phi^k_\eps$ of $\Delta_{g_\eps}$ as $\eps\to0$.
Each eigenvalue $\lambda^k$ satisfies the same estimate.
\end{proof}

\newcommand{\displacement}[1]{u^{#1}}
\newcommand{\velocity}[1]{v^{#1}}
\newcommand{\ii}{\mathrm{i}}
\newcommand{\dd}{\mathrm{d}}

\section{Original equations for gas giants}
\label{sec:gas-giant-physics}

Here, we present the extraction of the Laplace--Beltrami operator and acoustic wave operator from the system of equations describing the seismology on, and free oscillations of solar system gas giants. The system has been applied to studying the interiors of Saturn and Jupiter~\cite{DMFLX2021}. The original system is given explicitly for the displacement and contains implicitly the pressure; most of the work to extract the acoustic wave operator involves eliminating the displacement. Such an elimination appeared already in the study of inertial modes, that is, a reduction of the original system and invoking incompressibility leading to the Poincar\'{e} equation. Here, we follow the work of Prat, Ligni\`{e}res and Ballot~\cite{Prat_2016}.

\subsection{Acoustic-gravitational system of equations}


The displacement vector of a gas or liquid parcel between the unperturbed and perturbed flow is $\displacement{}$. The unperturbed values of pressure ($P$), density ($\rho$) and gravitational potential ($\Phi$) are denoted with a zero subscript.
The incremental Lagrangian stress formulation in the acoustic limit gives the equation of motion
\begin{multline}
   \rho_0 \partial_t^2 \displacement{}
     + 2 \rho_0 \Omega \times \partial_t \displacement{}
   = \nabla({\kappa \nabla \cdot {\displacement{}}})
   - \nabla (\rho_0 {\displacement{} \cdot \nabla (\Phi_0 + \Psi^s)})
\\
   + (\nabla \cdot (\rho_{0} \displacement{}))
     \nabla (\Phi_0 + \Psi^s)
   - \rho_0 \nabla\Phi' ,
\label{eq: MomentumConservation6}
\end{multline}
where the perturbed gravitational potential, $\Phi'$, solves
\begin{equation}
   \nabla^2 \Phi' = -4\pi G \nabla \cdot (\rho_0 \displacement{})
\label{eq: PoissonEqPert-2}
\end{equation}
and $\Psi^s$ denotes the centrifugal potential,
\begin{equation}
   \Psi^s = -\tfrac{1}{2} (\Omega^2 x^2 - (\Omega \cdot x)^2)
\label{eq: RotPot}
\end{equation}
($|\Omega|$ signifying the rotation rate of the planet). We may introduce the solution operator, $S$, for~\eqref{eq: PoissonEqPert-2} such that
\begin{equation}
   \Phi' = S(\rho_0 \displacement{}) .
\label{eq: PerturbGravOperator}
\end{equation}
We will use the shorthand notation,
\begin{equation}
   g_0' = -\nabla (\Phi_0 + \Psi^s).
\label{eq: PerturbGrav}
\end{equation}
A spherically symmetric manifold requires $\Omega = 0$ from well-posedness arguments.

\subsection{Brunt--V\"ais\"al\"a frequency}

We rewrite the first two terms on the right-hand side of~\eqref{eq: MomentumConservation6},
\begin{equation} \label{eq:rewrite-press}
   \nabla({\kappa \nabla \cdot {\displacement{}}})
   - \nabla (\rho_0 {\displacement{} \cdot \nabla (\Phi_0 + \Psi^s)})
   = \nabla [\kappa \rho_0^{-1} \,
             (\nabla \cdot (\rho_0 \displacement{})
     - \tilde{s} \cdot \displacement{})] ,
\end{equation}
in which
\begin{equation}
   \tilde{s} = \nabla \rho_0 - g_0' \frac{(\rho_0)^2}{\kappa} ,\quad
   \kappa = P_0 \gamma ;
\end{equation}
$\tilde{s}$ is related to the Brunt--V\"ais\"al\"a frequency, $N^2$,
\begin{equation}
   N^2 = \rho_0^{-1} (\tilde{s} \cdot g_0').
\label{eq: N2stilde}
\end{equation}
In~\eqref{eq:rewrite-press}, $\kappa \rho_0^{-1} \, (\nabla \cdot (\rho_0 \displacement{}) - \tilde{s} \cdot \displacement{})$ can be identified with the dynamic pressure, $-P$ say.
We recognize the acoustic wave speed,
\begin{equation}
   c^2 = \kappa \rho_0^{-1} .
\label{eq: SoundSpeed}
\end{equation}
Thus~\eqref{eq: MomentumConservation6} takes the form
\begin{multline}
   \partial_t^2 (\rho_0 \displacement{})
     + 2 \Omega \times \partial_t (\rho_0 \displacement{})
   = \nabla [c^2 \,
             (\nabla \cdot (\rho_0 \displacement{})
     - \rho_0^{-1} \tilde{s} \cdot (\rho_0 \displacement{}))]
\\
   + (\nabla \cdot (\rho_{0} \displacement{})) g_0'
     - \rho_0 \nabla\Phi' .
\label{eq: MCts}
\end{multline}
In~\eqref{eq: MCts} we can substitute~\eqref{eq: PerturbGravOperator} to arrive at an equation for $u$ containing a nonlocal contribution.

\subsection{Equivalent system of equations and Cowling approximation}

Writing $\velocity{} = \partial_t \displacement{}$ for the velocity, we obtain the following equivalent system of equations,
\begin{eqnarray}
   \partial_t \rho + \nabla \cdot (\rho_0 \velocity{}) &=& 0 ,
\\
   \partial_t (\rho_0 \velocity{})
     + 2 \Omega \times (\rho_0 \velocity{}) &=& -\nabla P
         + \rho g_0' - \rho_0 \nabla\Phi' ,
\\
   \partial_t P + \velocity{} \cdot \nabla P_0 &=& c^2 (\partial_t \rho         
     + \velocity{} \cdot \nabla \rho_0) ,
\end{eqnarray}
using that $\nabla P_0 = -\rho_0 g_0'$.
Well-posedness of the system of equations implies that
\begin{equation}
   \rho_0^{-1} (\tilde{s} \cdot (\rho_0 \displacement{})) g_0'
   = \frac{\rho_0^{-1} \tilde{s} \cdot g_0'}{|g_0'|^2}
     (g_0' \cdot (\rho_0 \displacement{})) .
\end{equation}
Upon inserting~\eqref{eq: N2stilde}, the third equation takes the form
\begin{equation}
   \partial_t P = c^2 \left(\partial_t \rho         
     + \frac{N^2}{|g_0'|^2} (g_0' \cdot (\rho_0 \velocity{})
       \right) .
\end{equation}
In the Cowling approximation the term $- \rho_0 \nabla\Phi'$ is dropped and the system reduces to
\begin{eqnarray}
   \partial_t \rho + \nabla \cdot (\rho_0 \velocity{}) &=& 0 ,
\label{eq:1}
\\[0.25cm]
   \partial_t (\rho_0 \velocity{})
     + 2 \Omega \times (\rho_0 \velocity{}) &=& -\nabla P
         + \rho g_0' ,
\label{eq:2}
\\
   \partial_t P &=& c^2 \left(\partial_t \rho         
     + \frac{N^2}{|g_0'|^2} (g_0' \cdot (\rho_0 \velocity{})
       \right) .
\label{eq:3}
\end{eqnarray}
This is identical to the system appearing in Prat \textit{et al.}~\cite{Prat_2016}.

\subsection{``Truncation'': Consistent boundary condition}

The free-surface boundary condition is given by the vanishing of the 
dynamic pressure (perturbation).
If $\rho_0$ and $c$ would not vanish at the boundary, we thus get the boundary condition
\begin{equation} \label{BC}
   (\kappa \nabla \cdot u + \rho_0 g_0' \cdot u)|_{\partial M_{\eps}} = 0 .
\end{equation}
(For comparison, the first term corresponds with the Lagrangian pressure perturbation.) This corresponds to taking the boundary condition slightly below the boundary rather than exactly at it.
We prove in proposition~\ref{Prop:truncation} that if the gas giant manifold is truncated just before the boundary, then the eigenvalues on this slightly smaller manifold $M_\eps$ converge to those of the true manifold $M$ at a specific rate.
This truncation has been widely used in computations~\cite{DMFLX2021}.

We have already noted the possibility of imposing certain types of boundary conditions when $\alpha < 2/n$, and this one here falls neatly into that framework. In particular, we can prove, just as in Section 4, that the domain of the operator augmented by this boundary condition is compactly contained in $L^2$, so that its spectrum is discrete. Furthermore, the spectra of the associated truncated problems converge at an estimable rate to the spectrum of this degenerate operator.

\subsection{Equation for the pressure and geometry}

We introduce
\[
   \nabla_z = \hat{\Omega} \cdot \nabla ,\quad
   \nabla_{\parallel} = (-\hat{g}_0') \cdot \nabla ,
\]
with unit vectors
\[
   \hat{\Omega} = \frac{\Omega}{|\Omega|} ,\quad
   \hat{g}_0' = \frac{g_0'}{|g_0'|}
\]
and
\[
   \nabla_{\perp} = \nabla + \hat{g}_0' \nabla_{\parallel} ,\quad
   \Delta_{\perp} = \nabla \cdot (\nabla_{\perp}) .
\]
Furthermore, $e_{\phi}$ is the unit vector in the direction of $\Omega \times g_0'$, noting that $\hat{\Omega}, \hat{g}_0', e_{\phi}$ form a non-orthogonal basis.

\begin{lemma}[{\cite{Prat_2016}}] \label{lem:Prat}
The time-Fourier-transformed pressure, $\hat{P}$, satisfies the equation  
\begin{multline} \label{eq:Lemma-Prat}
   \Delta \hat{P}
   - \frac{4}{\tau^2} \,
     \Omega \cdot \nabla (\Omega \cdot \nabla \hat{P})
\\
   - \frac{N^2}{(\tau^2 - 4 |\Omega|^2)}
     \Bigg[ \Delta \hat{P}
     - \frac{1}{|g_0'|^2} \, g_0' \cdot
     \nabla (g_0' \cdot \nabla \hat{P})
     - \frac{4}{\tau^2}
     \Omega \cdot \nabla (\Omega \cdot \nabla \hat{P})
\\
   - \frac{4 (\Omega \cdot g_0')^2}{\tau^2 |g_0'|^2} \,
     \Delta \hat{P}
   + \frac{4 (\Omega \cdot g_0')}{\tau^2 |g_0'|^2} \,
     (\Omega \cdot \nabla (g_0' \cdot \nabla \hat{P})
       + g_0' \cdot \nabla (\Omega \cdot \nabla \hat{P}))
   \Bigg]
\\
   + \frac{1}{\tau^2 (\tau^2 - 4 |\Omega|^2)}
     \mathscr{V} \cdot \nabla \hat{P}    
   + \frac{1}{c^2 \tau^4 (\tau^2 - 4 |\Omega|^2)} \,
     \mathscr{R}
\\
   + \frac{1}{c^2} \Bigg[\tau^2 \Bigg(1 - \frac{4}{\tau^2} |\Omega|^2\Bigg)
     + c^2 \widehat{\mathcal{M}} \,
     \nabla \cdot \Bigg(
     \frac{(4 \tau^{-2} (g_0' \cdot \Omega) \Omega
      - g_0')}{c^2 \widehat{\mathcal{M}}}\Bigg)\Bigg]
     \widehat{P} = 0 ,
\end{multline}
where $\mathscr{R}$ and $\mathscr{V}$ are given below. In coordinates relative to the above mentioned non-orthogonal basis, the terms with leading, second-order spatial derivatives take the form $\Delta \hat{P} - \tau^{-2} (4 |\Omega|^2 \, \nabla_z^2 \hat{P} + N^2 \Delta_{\perp} \hat{P})$; the leading, second-order term in $\tau$ is given by $c^{-2} \tau^2 \hat{P}$. Thus one identifies, to leading order, the acoustic wave operator on the one hand and an equation like Poincar\'{e}'s equation in the (axi)symmetric case~\cite{RieutordNoui_1999} on the other hand.
\end{lemma}

\noindent
For clarity, we summarize the proof of this lemma. To eliminate $\velocity{}$ from the system of equations, one takes $2 \Omega \times$ and $2 \Omega \, \cdot$ of~\eqref{eq:2} and applies a time derivative to the resulting equations. Upon taking another time derivative, and substituting the second resulting equation in the first, one obtains the equation
\begin{multline}
   \mathcal{L}(\rho_0 \velocity{}) = -4 (\Omega \cdot \nabla P) \, \Omega
     + 4 \rho (\Omega \cdot g_0') \, \Omega
\\
     - \nabla \partial_t^2 P + (\partial_t^2 \rho) g_0'
     + 2 \Omega \times \nabla \partial_t P
     - 2 (\partial_t \rho) \Omega \times g_0' ,
\label{eq:P-1}
\end{multline}
where
\[
   \mathcal{L} = \partial_t^3 + 2 |\Omega|^2 \partial_t .
\]
With this operator, \eqref{eq:1} implies
\begin{equation}
   \partial_t \mathcal{L}(\rho) =
        - \nabla \cdot \mathcal{L}(\rho_0 \velocity{})
\label{eq:P-2}
\end{equation}
and~\eqref{eq:3} implies
\begin{equation}
   \partial_t \mathcal{L}(P) = c^2 \partial_t \mathcal{L}(\rho)
     + \beta g_0' \cdot \mathcal{L}(\rho_0 \velocity{}) ,
\label{eq:P-3}
\end{equation}
where
\[
   \beta = \frac{c^2 N^2}{|g_0'|^2} .
\]
Substituting~\eqref{eq:P-1}) into~\eqref{eq:P-3} gives
\begin{multline}
   \partial_t \mathcal{L}(P)
   + \beta [4 (\Omega \cdot \nabla P) (g_0' \cdot \Omega) 
     + g_0' \cdot \nabla \partial_t^2 P
\\
     + 2 (\Omega \times g_0') \cdot \nabla \partial_t P]
   = c^2 \partial_t \mathcal{L}(\rho)
     + \beta[\rho (2 \Omega \cdot g_0')^2
       + (\partial_t^2 \rho) |g_0'|^2] .
\label{eq:P-4}
\end{multline}
Using the definition of $\mathcal{L}$, one may extend the operator notation to $c^2 \mathcal{M}(\rho)$ for the right-hand side of this equation, with
\[
   \mathcal{M} = \partial_t^4 + (4 |\Omega|^2 + N^2) \partial_t^2
       + \frac{4 N^2 (\Omega \cdot g_0')^2}{|g_0'|^2} .
\]
Introducing the dual, $\tau$, of $\ii \partial_t$, one writes
\[
   \widehat{\partial_t \mathcal{L}} = \tau^2 (\tau^2 - 4 |\Omega|^2),
\quad
   \widehat{\mathcal{M}} =  \tau^4 - (4 |\Omega|^2 + N^2) \tau^2
       + \frac{4 N^2 (\Omega \cdot g_0')^2}{|g_0'|^2}
\]
for the relevant symbols, noting that $N^2$ and $g_0'$ are dependent on the coordinates. Equation~\eqref{eq:P-4} then gives
\begin{equation}
   \hat{\rho} = \frac{\widehat{\partial_t \mathcal{L}} \, \hat{P}
     + \beta [4 (g_0' \cdot \Omega) \, \Omega
     - \tau^2 g_0'
     - 2 \ii \tau (\Omega \times g_0')] \cdot \nabla \hat{P}}{
                c^2 \widehat{\mathcal{M}}} .
\label{eq:P-5}
\end{equation}
Taking the divergence of~\eqref{eq:P-1} yields
\begin{multline}
   \nabla \cdot \mathcal{L}(\rho_0 \hat{\velocity{}})
   = \tau^2 \Delta \hat{P}
\\
     - 4 \Omega \cdot \nabla (\Omega \cdot \nabla \hat{P})
     + \nabla \cdot [\hat{\rho} \, (
       4 (g_0' \cdot \Omega) \Omega
       - \tau^2 g_0'
       + 2 \ii \tau (\Omega \times g_0'))] .
\label{eq:P-6}
\end{multline}
Here, it was used that $\Omega$ is a constant vector (signifying uniform rotation). One then considers~\eqref{eq:P-2} and substitutes~\eqref{eq:P-6} to obtain an equation for $\hat{P}$ upon using~\eqref{eq:P-5} on the left-hand side:
\begin{multline}
   \tau^2 \Delta \hat{P}
   - 4 \Omega \cdot \nabla (\Omega \cdot \nabla \hat{P})
\\
   + \nabla \cdot \Bigg[
     \frac{(\widehat{\partial_t \mathcal{L}} \, \hat{P}
      + \beta [4 (g_0' \cdot \Omega) \Omega
      - \tau^2 g_0'
      - 2 \ii \tau (\Omega \times g_0')] \cdot \nabla \hat{P})}{
     c^2 \widehat{\mathcal{M}}} 
\\     
   \cdot (4 (g_0' \cdot \Omega) \Omega
      - \tau^2 g_0'
      + 2 \ii \tau (\Omega \times g_0')) \Bigg]
\\
   + \frac{\widehat{\partial_t \mathcal{L}} \,
     (\widehat{\partial_t \mathcal{L}} \, \hat{P}
      + \beta [4 (g_0' \cdot \Omega) \Omega
     - \tau^2 g_0'
     - 2 \ii \tau (\Omega \times g_0')] \cdot \nabla \hat{P})}{
     c^2 \widehat{\mathcal{M}}} = 0 .
\label{eq:P-7}
\end{multline}
As $g_0'$ derives from a potential (cf.~\eqref{eq: PerturbGrav}), it follows that
\[
   \nabla \cdot (\Omega \times g_0') = 0 .
\]
Then
\begin{multline}
   \tau^2 \Delta \hat{P}
   - 4 \Omega \cdot \nabla (\Omega \cdot \nabla \hat{P})
\\
   + \nabla \cdot \Bigg[
     \frac{(\widehat{\partial_t \mathcal{L}} \, \hat{P}
      + \beta [4 (g_0' \cdot \Omega) \Omega
      - \tau^2 g_0'] \cdot \nabla \hat{P})}{
     c^2 \widehat{\mathcal{M}}} \,
     (4 (g_0' \cdot \Omega) \Omega
      - \tau^2 g_0') \Bigg]
\\
   + \frac{\widehat{\partial_t \mathcal{L}} \,
     (\widehat{\partial_t \mathcal{L}} \, \hat{P}
      + \beta [4 (g_0' \cdot \Omega) \Omega
     - \tau^2 g_0'] \cdot \nabla \hat{P})}{
     c^2 \widehat{\mathcal{M}}}
     + \frac{1}{c^2 \widehat{\mathcal{M}}} \mathscr{R} = 0 ,
\label{eq:P-8}
\end{multline}
where
\begin{multline}
   \mathscr{R} = -2 \ii \tau c^2 \widehat{\mathcal{M}} \,
     \nabla \cdot \Bigg[
     \frac{\beta}{c^2 \widehat{\mathcal{M}}}
       ((\Omega \times g_0') \cdot \nabla \hat{P}) \,
     (4 (\Omega \cdot g_0') \, \Omega 
     + 2 \ii \tau \Omega \times g_0' - \tau^2 g_0') \Bigg]
\\
   + 2 \ii \tau c^2 (\Omega \times g_0') \cdot \nabla \Bigg[   
     \frac{\beta}{c^2 \widehat{\mathcal{M}}}
       (4 (\Omega \cdot g_0') \, \Omega - \tau^2 g_0')
         \cdot \nabla \hat{P} \Bigg] 
\\
   + 2 \ii \tau c^2 \widehat{\partial_t \mathcal{L}} \,
    (\Omega \times g_0') \cdot
     \Bigg[\nabla\Bigg(\frac{\hat{P}}{c^2 \widehat{\mathcal{M}}}\Bigg)
       - \frac{\beta}{c^2 \widehat{\mathcal{M}}} \nabla \hat{P}\Bigg]
\label{eq:P-9}
\end{multline}
represents the sum of terms containing $2 \ii \tau (\Omega \times g_0') \, \cdot$. It is noted that in the axisymmetric case,
\[
   (\Omega \times g_0') \cdot \nabla \Bigg(
   \frac{1}{c^2 \widehat{\mathcal{M}}} \Bigg) = 0 ,
\]
and that in a polytropic model (see Subsection~\ref{sec:physics-intro}), $\beta$ is a constant, which simplifies the computations. Equation~\eqref{eq:P-8} can be rewritten as
\begin{multline}
   \widehat{\mathcal{M}} \, \Delta \hat{P}
   - 4 (\tau^2 - (4 |\Omega|^2 + N^2)) \,
     \Omega \cdot \nabla (\Omega \cdot \nabla \hat{P})
\\[0.25cm]
   + \frac{\beta \tau^2}{c^2} \, g_0' \cdot
     \nabla (g_0' \cdot \nabla \hat{P})
   - \frac{4 \beta}{c^2} \, (\Omega \cdot g_0') \,
     [\Omega \cdot \nabla (g_0' \cdot \nabla \hat{P})
       + g_0' \cdot \nabla (\Omega \cdot \nabla \hat{P})]
\\
   + \frac{4}{c^2 \tau^2} \Bigg\{\Bigg[(1 + \beta)
     \widehat{\partial_t \mathcal{L}} + c^2 \widehat{\mathcal{M}} \,
     \nabla \cdot \Bigg(\frac{\beta (4 (g_0' \cdot \Omega) \Omega
       - \tau^2 g_0')}{c^2 \widehat{\mathcal{M}}}\Bigg)\Bigg] \,
         (g_0' \cdot \Omega)
\\
      + \frac{\beta}{2}
      \Omega \cdot \nabla (g_0' \cdot (4 (g_0' \cdot \Omega) \Omega
        - \tau^2 g_0')) \Bigg\} \,
      \Omega \cdot \nabla \hat{P}    
\\
   - \frac{1}{c^2} \Bigg\{(1 + \beta)
     \widehat{\partial_t \mathcal{L}} + c^2 \widehat{\mathcal{M}} \,
     \nabla \cdot \Bigg(\frac{\beta (4 (g_0' \cdot \Omega) \Omega
       - \tau^2 g_0')}{c^2 \widehat{\mathcal{M}}}\Bigg)\Bigg\} \,
         g_0' \cdot \nabla \hat{P}    
\\
   + \frac{\widehat{\partial_t \mathcal{L}}}{c^2 \tau^2} \Bigg[
     \widehat{\partial_t \mathcal{L}} + c^2 \widehat{\mathcal{M}} \,
     \nabla \cdot \Bigg(\frac{(4 (g_0' \cdot \Omega) \Omega
      - \tau^2 g_0')}{c^2 \widehat{\mathcal{M}}}\Bigg)\Bigg] \,
     \hat{P}
   + \frac{1}{c^2 \tau^2} \, \mathscr{R} = 0 .
\label{eq:P-11}
\end{multline}
The sum of the two terms containing factors in between braces allow the shorthand notation $\mathscr{V} \cdot \nabla \hat{P}$: 
\begin{multline}
   (\tau^2 - 4 |\Omega|^2) \, \Delta \hat{P}
   - (\tau^2 - 4 |\Omega|^2) \, \frac{4}{\tau^2} \,
     \Omega \cdot \nabla (\Omega \cdot \nabla \hat{P})
\\
   - N^2 \Bigg[ \Delta \hat{P}
   - \frac{1}{|g_0'|^2} \, g_0' \cdot
     \nabla (g_0' \cdot \nabla \hat{P})
   - \frac{4}{\tau^2}
     \Omega \cdot \nabla (\Omega \cdot \nabla \hat{P})
\\
   - \frac{4 (\Omega \cdot g_0')^2}{\tau^2 |g_0'|^2} \,
     \Delta \hat{P}
   + \frac{4 (\Omega \cdot g_0')}{\tau^2 |g_0'|^2} \,
     (\Omega \cdot \nabla (g_0' \cdot \nabla \hat{P})
       + g_0' \cdot \nabla (\Omega \cdot \nabla \hat{P}))
   \Bigg]
   + \frac{1}{\tau^2} \mathscr{V} \cdot \nabla \hat{P}    
\\
   + \frac{\widehat{\partial_t \mathcal{L}}}{c^2 \tau^4} \Bigg[
     \widehat{\partial_t \mathcal{L}} + c^2 \tau^2 \widehat{\mathcal{M}} \,
     \nabla \cdot \Bigg(
       \frac{(4 \tau^{-2} (g_0' \cdot \Omega) \Omega
       - g_0')}{c^2 \widehat{\mathcal{M}}}\Bigg)\Bigg] \,
     \hat{P}
   + \frac{1}{c^2 \tau^4} \, \mathscr{R} = 0 .
\label{eq:P-12}
\end{multline}
One then divides the equation by $(\tau^2 - 4 |\Omega|^2)$; the factor in front of $\hat{P}$ then takes the form
\begin{equation}
   \frac{1}{c^2} \Bigg[\tau^2 \Bigg(1 - \frac{4}{\tau^2} |\Omega|^2\Bigg)
     + c^2 \widehat{\mathcal{M}} \,
     \nabla \cdot \Bigg(
     \frac{(4 \tau^{-2} (g_0' \cdot \Omega) \Omega
      - g_0')}{c^2 \widehat{\mathcal{M}}}\Bigg)\Bigg]
    .
\end{equation}
This results in equation~\eqref{eq:Lemma-Prat}.

\subsection{Propagation of singularities}

The propagation of singularities depends only on the leading order part of the system of equations. Ignoring lower order terms, equation~\eqref{eq: MCts} reads
\begin{equation}
   \partial_t^2 (\rho_0 \displacement{})
   - \nabla [c^2  \nabla \cdot (\rho_0 \displacement{})]
   =
   0
\end{equation}
and the principal symbol at $(t,x;\tau,\xi)$ is $\tau^2\mathrm{Id}-c^2(x)\xi\xi^T$.
From the way the matrix $\xi\xi^T$ acts we may read that pressure singularities propagate but shear ones do not. Pressure waves (``polarized'' along the momentum $\xi$) follow the geodesics of the isotropic sound speed $c$ just as the solutions of the scalar wave equation for pressure $(\partial_t^2 - c^2 \Delta) P = 0$ as extracted from the original system in Lemma~\ref{lem:Prat}. Therefore, if only the travel times of singularities are concerned, it suffices to model a gas planet with a scalar wave equation.

The parametrix construction outlined in Section~\ref{sec:4} is a very flexible one. Although we have used it to analyze the simpler operator $\Delta_g$ studied
in the rest of this paper, the operator appearing in ~\eqref{eq:Lemma-Prat} is a perturbation of such a Laplacian, for appropriately defined gas-giant metric $g$,
with all extra terms being of lower order in the sense of this calculus of degenerate operators.  In other words, it is possible, just as easily, to construct a
parametrix for this operator in the $0$-pseudodifferential calculus, and to derive the same sorts of conclusions as we have discussed for the Laplacian.
Furthermore, the lower order terms here are all compact relative to the main part of this operator, hence do not affect the discreteness of the spectrum, 
but do cause the usual sorts of perturbations to the spectrum caused by any such compact perturbations.


\bibliographystyle{alpha}
\bibliography{references}

\end{document}